\documentclass{article}
\usepackage{graphicx}
\usepackage{authblk}
\usepackage{hyperref}
\hypersetup{
	pdftitle={Optimization of piecewise smooth shapes under uncertainty using the example of Navier--Stokes flow},
	pdfauthor={Caroline Geiersbach, Tim Suchan, and Kathrin Welker}
}
\usepackage[T1]{fontenc}
\usepackage[utf8]{inputenc}
\usepackage{lmodern} 

\providecommand{\keywords}[1]
{
	\small	
	\textbf{\textit{Keywords. }} #1
}
\providecommand{\ams}[1]
{
	\small	
	\textbf{\textit{AMS subject classifications. }} #1
}

\usepackage{amsmath,amssymb,amscd,latexsym} 
\usepackage{upgreek}
\usepackage{bbm} 
\usepackage{amsfonts}
\usepackage{mathrsfs}
\usepackage{cancel}
\usepackage{mathtools} 
\usepackage{stmaryrd} 
\usepackage{afterpage} 
\usepackage[numbers, sort]{natbib}

\usepackage{algorithm}
\usepackage[noend]{algpseudocode}
\algnewcommand{\IfThenElse}[3]{
	\State \algorithmicif\ #1\ \algorithmicthen\ #2\ \algorithmicelse\ #3}

\usepackage{color}
\usepackage{graphicx}
\graphicspath{{figures/}{tikzs/}}
\usepackage{overpic}
\usepackage{subcaption}
\usepackage{placeins} 
\usepackage{footmisc}  
\usepackage{tablefootnote}
\newlength\figureheight 
\newlength\figurewidth
\graphicspath{{figures/}{tikzs/}{svgs/}}
\usepackage{svg}
\svgpath{{svgs/}}
\usepackage{tikz}
\usetikzlibrary{shapes.geometric,calc,positioning}
\usepackage{tikz-cd}
\usepackage{tikz-3dplot}
\usepackage{pgfplots}				
\usepackage{varwidth}
\newlength\mysvgwidth
\newlength\boxwidth
\newlength\svgwidth
\pgfplotsset{compat=newest,grid style={dotted,black}}
\usepgfplotslibrary{external} 			
\usepgfplotslibrary{patchplots}
\ifpdf
  \DeclareGraphicsExtensions{.eps,.pdf,.png,.jpg}
\else
  \DeclareGraphicsExtensions{.eps}
\fi

\newcommand{\R}{{\mathbb{R}}}

\newcommand{\E}{{\mathbb{E}}}
\newcommand{\pP}{{\mathbb{P}}}

\renewcommand{\vec}[1]{\boldsymbol{#1}}
\newcommand{\Deriv}{\mathrm{d}_{\hat{\vu}}}
\newcommand{\partDeriv}{\mathrm{d}_{i}}
\newcommand{\partMatDeriv}{\mathrm{d}_{m_i}}
\newcommand{\matDeriv}{\mathrm{d}_{m}}

\DeclareMathOperator{\Divv}{div}
\newcommand{\Div}[1]{\Divv{(#1)}}

\DeclareMathOperator{\vol}{vol}
\DeclareMathOperator{\bary}{bary}

\renewcommand{\d}{\,\mathrm{d}}
\newcommand{\dx}{\,\mathrm{d} \vec{x}}
\newcommand{\ds}{\,\mathrm{d} \vec{s}}
\newcommand{\vx}{\vec{x}}

\newcommand{\velocity}{\vec{v}}
\newcommand{\pressure}{p}

\newcommand{\adjvelocity}{\vec{\varphi}}
\newcommand{\adjpressure}{\psi}

\newcommand{\vn}{\vec{n}}

\newcommand{\vu}{\vec{u}}

\newcommand{\cU}{\mathcal{U}}

\newcommand{\GD}{\Gamma_D}
\newcommand{\GN}{\Gamma_N}
\newcommand{\M}{M_s(\mathcal{U}^N)}

\newcommand{\intD}[1]{\int_{D_{\hat{\vu}}} #1 \dx}
\newcommand{\intDeltai}[1]{\int_{\Delta_i} #1 \dx}
\newcommand{\intDi}[1]{\int_{D_{\hat{u}_i}} #1 \dx}
\newcommand{\intui}[1]{\int_{\hat{u}_i} #1 \ds}

\newcommand{\restr}[2]{#1\rvert_{#2}}

\newtheorem{definition}{Definition}[section]

\newtheorem{remark}{Remark}

\begin{document}


\title{Optimization of piecewise smooth shapes under uncertainty using the example of Navier--Stokes flow}


\author{Caroline Geiersbach}
\affil{Weierstrass Institute, Mohrenstraße 39, 10117 Berlin, Germany
\href{mailto:caroline.geiersbach@wias-berlin.de}{\ttfamily caroline.geiersbach@wias-berlin.de}\vspace*{6pt}}

\author{Tim Suchan}
\affil{Helmut-Schmidt-Universität / Universität der Bundeswehr Hamburg, Holstenhofweg 85, 22043 Hamburg, Germany, \href{mailto:suchan@hsu-hh.de}{\ttfamily suchan@hsu-hh.de}\vspace*{6pt}}

\author{Kathrin Welker}
\affil{Technische Universität Bergakademie Freiberg, Akademiestraße 6, 09599 Freiberg, Germany, \href{mailto:Kathrin.Welker@math.tu-freiberg.de}{\ttfamily Kathrin.Welker@math.tu-freiberg.de}}

\date{\today}

\maketitle

\begin{abstract}
We investigate a complex system involving multiple shapes to be optimized in a domain, taking into account geometric constraints on the shapes and uncertainty appearing in the physics. We connect the differential geometry of product shape manifolds with multi-shape calculus, which provides a novel framework for the handling of piecewise smooth shapes. This multi-shape calculus is applied to a shape optimization problem where shapes serve as obstacles in a system governed by steady state incompressible Navier--Stokes flow. Numerical experiments use our recently developed stochastic augmented Lagrangian method and we investigate the choice of algorithmic parameters using the example of this application. 
\end{abstract}

\newline
\ams{
	49Q12, 49K20, 53C15
}

\section{Introduction} \label{sec:intro}

Shape optimization is concerned with identifying shapes, or subsets of $\R^d$, behaving in an optimal way with respect to a given physical system. Many problems of interest involve a system in the form of a partial differential equation (PDE), the solution of which depends on  one or more shapes defining the domain. Some applications involve additional geometric constraints on the shapes, leading to nonsmooth problems that cannot be solved using standard descent-type algorithms. Moreover, there has been an increasing interest in incorporating uncertain parameters or inputs in shape optimization models. 

Shape optimization is commonly applied  in engineering in order to optimize shapes.
Theory and algorithms in shape optimization can be based on techniques from differential geometry, e.g., a Riemannian manifold structure can be used to define the distances of two shapes.
Thus, shape spaces are of particular interest in shape optimization.
The shape space
$B_e(S^1, \mathbb{R}^2)$
briefly investigated in~\cite{MichorMumford1} is an important example of a smooth manifold allowing a Riemannian structure, where the term smooth shall refer to infinite differentiability in this paper. This shape space is considered in recent publications (cf., e.g., \cite{Geiersbach2021,geiersbach2023stochastic,Schulz,Schulz2016,Schulz2016a}),  but  $B_e(S^1, \mathbb{R}^2)$ is in general not sufficient to carry out optimization algorithms on piecewise smooth shapes. In particular, a piecewise smooth shape is often encountered as an optimal shape for fluid-mechanical problems, see e.g.~\cite{Pironneau1973}.
Some effort has been put into constructing a shape space that contains non-smooth shapes. A shape space which is a diffeological space is defined in \cite{Welker2021}.
Recently, a space containing shapes in $\R^2$ that can be identified with a Riemannian product manifold but at the same time admits piecewise smooth curves as elements was constructed \cite{Pryymak2023}. In this paper, we focus on this shape space.
In many applications in shape optimization, it is also desirable to consider multiple shapes to be optimized \cite{Albuquerque2020,Cheney1999,Kwon2002}. A first approach for optimizing multiple smooth shapes was presented in~\cite{Geiersbach2022} and applied in \cite{geiersbach2023stochastic,Pryymak2023}.

An early work in incorporating uncertainty in shape optimization \cite{Conti2008} placed stochastic models in the context of stochastic programming, a classical topic concerned with optimization subject to uncertainty. A helpful guide on models and methods in shape optimization under uncertainty is \cite{Martinez-Frutos2016}. For problems with a finite or low stochastic dimension, or in the special case where the problem's structure can be exploited \cite{Dambrine2015}, the stochastic quantity may be exactly represented or discretized using the stochastic Galerkin method or polynomial chaos; cf.~\cite{Atwal2012,Conti2018,Schillings2011}. For larger dimensions, ensemble-based approaches are typically needed. This includes stochastic collocation \cite{Martinez-Frutos2016} and Monte Carlo-based methods. In the context of optimization, if a random sample is generated and the original problem is replaced by a proxy problem with this sample, one speaks of sample average approximation. Another approach is stochastic approximation, which involves dynamically sampling as part of the underlying method. Stochastic approximation for shape spaces was proposed in \cite{Geiersbach2021} and further developed for multi-shape problems in \cite{Geiersbach2022}. 

In this paper, we handle a stochastic shape optimization problem of the form 
\begin{equation}
\label{eq:SO-problem-abstract-extended}
\begin{aligned}
&\min_{\vec{u} \in \mathcal{M}} \, \left \lbrace \E[J(\vec{u},\vec{\xi})] = \int_\Omega J(\vec{u},\vec{\xi}(\omega)) \d \pP(\omega)  \right\rbrace \\
&\text{subject to (s.t.)} \quad h_i(\vec{u}) \leq 0 \quad \forall i \in \{1, \dots, n\}. 
\end{aligned}
\end{equation}
Here, $\vu \coloneqq \left( u_1, \ldots, u_s \right)$ is a vector of $s$  shapes in a shape space $\mathcal{M}$, $\vec{\xi}$ is a random vector with probability measure $\pP\colon \Omega \rightarrow [0,1]$, and the functions $h_i \colon \mathcal{M} \rightarrow \R$ are constraints on the shapes. We focus on further developing shape calculus to handle problems involving the novel manifold of piecewise smooth shapes from \cite{Pryymak2023}. 
We define the multi-material derivative as well as the multi-shape derivative and its stochastic version and connect it with this manifold. Furthermore, we compute the multi-shape derivative for a system subject to Navier--Stokes flow under uncertainty with deterministic geometric constraints on the shapes. To solve this problem computationally, we use the stochastic augmented Lagrangian method recently proposed in \cite{geiersbach2023stochastic}. The method uses the framework of stochastic approximation on shape spaces from \cite{Geiersbach2021,Geiersbach2022}, 
combined with an augmented Lagrangian procedure based on \cite{Steck2018,Kanzow2018} for updating penalty parameters and Lagrange multipliers. In contrast to \cite{geiersbach2023stochastic}, we do not require infinitely smooth shapes. It is, however, still possible to define a stochastic multi-shape derivative for nonsmooth shapes, which is used in the inner loop procedure of the stochastic augmented Lagrangian method. 

The paper is structured as follows. In Section~\ref{sec:Recap}, we present the shape space of piecewise smooth shapes and connect it with multi-shape calculus. In particular, we define the (stochastic) multi-shape derivative, which is needed for the application. In Section~\ref{sec:specificProblem}, we describe an example multi-shape optimization problem and compute its stochastic multi-shape derivative using a novel multi-material derivative approach. In Section~\ref{sec:numericalResults}, we detail the stochastic augmented Lagrangian method for multi-shape optimization and present numerical results demonstrating the convergence rates from our theory and compare results from stochastic and deterministic optimization. We conclude with a short discussion in Section~\ref{sec:conclusions}.

\section{Multi-shape calculus for optimizing piecewise smooth shapes}
\label{sec:Recap} 
In  Section~\ref{sec:PSManifold}, we introduce a space containing shapes in $\R^2$ that can be identified with a Riemannian product manifold but at the same time admits piecewise smooth curves as elements. This shape space was first defined in \cite{Pryymak2023} and is suitable for our application in  Section~\ref{sec:specificProblem}, a fluid-mechanical problem constrained by the Navier--Stokes equations. In Section~\ref{subsec:multi-shape-derivative}, we define the (stochastic) multi-shape derivative for this shape space.

\subsection{The product shape manifold of piecewise smooth shapes}
\label{sec:PSManifold}
Let $\mathcal{U}_1, \dots, \mathcal{U}_N$  be manifolds equipped with Riemannian metrics $G^{i'} = (G_u^{i'})_{u \in \mathcal{U}_{i'}}$, $i'=1,\dots,N$. We define the Riemannian product manifold $\mathcal{U}^N$ by
$  \mathcal{U}^N \coloneqq \mathcal{U}_1 \times \cdots \times \mathcal{U}_N= \prod_{i'=1}^{N}\mathcal{U}_{i'}$.
The product metric $\mathcal{G}^N$ to the corresponding product shape space $\mathcal{U}^N$ can be defined via 
$\mathcal{G}^N=\sum_{i'=1}^N \pi_{i'}^\ast \mathcal{G}^{i'}$,
where $\pi_{i'}^{*}$ are the pullbacks associated with canonical projections $\pi_{i'}\colon \cU^N\to \cU_{i'}$, $i'=1, \dots, N$ (cf. \cite{Geiersbach2022}). 
As in \cite{Pryymak2023},  we define the $s$-dimensional shape space on $\mathcal{U}^N$ by 
\begin{align*}
M_s(\mathcal{U}^N)\coloneqq \Big\lbrace \vec{u}=(u_1, \dots ,u_s) \mid \,& u_i \in  \prod\limits_{l=k_i}^{k_i+n_i-1} \mathcal{U}_{l} \,\, \forall i=1,\ldots,s,\, \sum\limits_{i=1}^s n_i = N \text{ and }\\
&k_1=1, \, k_{i+1}=k_i+n_i
\,\forall i =1,\dots,s-1  \Big\rbrace.
\end{align*}
An element in $M_s(\mathcal{U}^N)$ is defined as a vector of $s$ shapes $u_1,\dots,u_s$, where each shape $u_i$ is an element of the product of $n_i$ smooth manifolds.
Note that any element $\vec{u}=(u_1,\dots,u_s)\in M_s(\mathcal{U}^N)$ can be understood as an element $\tilde{\vec{u}}=(\tilde{u}_1,\dots, \tilde{u}_N)\in\mathcal{U}^N$. Thus, we identify
$T_{\vec{u}} M_s(\mathcal{U}^N) = T_{\tilde{\vec{u}}} \mathcal{U}^N$  and the Riemannian metric $\mathcal{G}=(\mathcal{G}_{\vu})_{\vu \in \M}$ with
$\mathcal{G}^N = (\mathcal{G}^N_{\tilde{\vec{u}}})_{\tilde{\vu} \in \mathcal{U}^N}.$

Since we are interested in optimizing piecewise smooth shapes, we restrict the choice of shapes in $M_s(\mathcal{U}^N)$ to piecewise smooth shapes that are glued together. More precisely, we assume that each shape $(u_1,\dots,u_s)\in M_s(\mathcal{U}^N)$ is single closed: either $u_i \in  B_e(S^1, \mathbb{R}^2)$ or $u_i =(u_{k_i},\dots,u_{k_i+n_i-1}) \in   (B_e([0,1], \mathbb{R}^2))^{n_i}$  with the additional conditions 
\begin{equation}
\label{eq:glue-condition}
\begin{split}
& u_i\colon [0,1)\to \mathbb{R}^2 \text{ injective with }\\
& u_{k_i+h}(1)=u_{k_i+h+1}(0) \, \forall h = 0,\dots,n_i-2 \text{ and } u_{k_i}(0)=u_{k_i+n_i-1}(1).
\end{split}
\end{equation}
Here, the shape spaces are defined by
\begin{align*}
B_e(S^1, \mathbb{R}^2) &\coloneqq \text{Emb}(S^1, \mathbb{R}^2) /\text{Diff}(S^1),\\
B_e([0,1], \mathbb{R}^2) &\coloneqq\text{Emb}([0,1], \mathbb{R}^2) /\text{Diff}([0,1]),
\end{align*}
where $\text{Emb}(\cdot, \mathbb{R}^2)$ denotes the set of all embeddings into $\R^2$ and $\text{Diff}(\cdot)$ is the set of all diffeomorphisms.
Since a smooth curve\footnote{Throughout this paper, a curve $c$ is to be understood as a continuous function $c\colon I \rightarrow \R^2$, where $I \neq \emptyset$ is an interval in $\R$.} can be defined by an embedding (cf., e.g., \cite[Chapter 2, Definition 2.22]{Kuehnel}), these spaces represent all simple closed smooth curves in $\mathbb{R}^2$ and all simple open smooth curves in $\mathbb{R}^2$, respectively.
Both spaces are smooth manifolds allowing Riemannian structures (cf., e.g., \cite{MichorMumford1,Michor1980}).
A closed curve with kinks is interpreted as a glued-together curve of open smooth curves, i.e.,  elements of $B_e([0,1], \mathbb{R}^2)$.
In contrast, a shape defined as an element in $B_e(S^1, \mathbb{R}^2)$ has no kinks in the shape geometry; in the following application, we will not consider this case. We focus on glued-together piecewise smooth shapes. For this, it will be convenient to define 
\begin{align*}
M_s^c &\coloneqq \Big\lbrace \vec{u} \in M_s(B_e([0,1], \mathbb{R}^2)^N) \mid u_i =(u_{k_i},\dots,u_{k_i+n_i-1}) \in   (B_e([0,1], \mathbb{R}^2))^{n_i}   \\
&\quad \quad \quad  \text{ with } \eqref{eq:glue-condition} \, \forall i=1,\dots,s,\, \sum_{i=1}^s n_i=N \Big\rbrace.
\end{align*}

\begin{remark}\label{rem:KinksInShapes}
An element of $B_e(S^1, \R^2)$, by definition, is smooth and so does not contain kinks. The product shape manifold described above from \cite{Pryymak2023} allows for a larger set of possible shapes and also includes shapes from $B_e(S^1, \R^2)$. In computations, the number $n_i$ of glued-together smooth shapes, i.e., the maximum number of kinks of each single closed shape $u_i$ with $i=1,\ldots,s$, can be chosen to be the number of nodes belonging to the discretization of this shape. In this way, elements cannot leave the shape space over the course of optimization.
\end{remark}


\subsection{Multi-shape derivative}
\label{subsec:multi-shape-derivative}
As already mentioned, we are interested in optimizing with respect to multiple glued-together piecewise smooth shapes. For this, we will need to apply the definition of the partial shape derivative and the multi-shape derivative given in \cite{Geiersbach2022} to our setting. Let $\vec{u} = (u_1, \dots, u_s) \in M_s^c$. We identify a curve with its image; for each glued-together piecewise smooth shape $u_i$, we define the corresponding subset in $\R^2$ by \
\begin{equation}
\label{eq:subset-for-shape}
\begin{split}
\hat{u}_i \coloneqq \{ & \vec{x} \in \{u_{k_i}([0,1]),\dots,u_{k_i+n_i-1}([0,1])\} \mid \\ & (u_{k_i},\dots,u_{k_i+n_i-1}) \in   (B_e([0,1], \mathbb{R}^2))^{n_i} \text{ with } \eqref{eq:glue-condition} \}
\end{split}
\end{equation}
and the set of the corresponding vector of subsets of $\R^2$ by 
\begin{align*}
\hat{M}_s^c &= \{ \hat{\vu}=(\hat{u}_1, \dots, \hat{u}_s) \subset (\R^2)^s \}.
\end{align*} 
 We define a function $\mathbb{U} \colon M^c_s \rightarrow \hat{M}^c_s, \vu \mapsto \hat{\vu}$ corresponding to the identification described in \eqref{eq:subset-for-shape}. In this way, we can identify a function $H\colon M_s^c \rightarrow \R$ defined on the manifold with the function $\hat{H}\colon \hat{M}_s^c \rightarrow \R$ defined on subsets of $\R^2$ via $H(\vu) = \hat{H}(\mathbb{U}(\vec{u}))$ for all $\vu \in M_s^c$.

For the definition of the partial shape derivative, we need transformations of the individual subsets $\hat{u}_i$. Let $D \subset \R^2$ be a domain containing $s$ non-overlapping  shapes $\hat{\vu} \in \hat{M}_s^c$. Let $\Delta_1,\ldots,\Delta_s$ be a partition of $D$ into $s$
non-empty, connected, bounded sets with Lipschitz boundaries
such that $\hat{u}_i \subset \Delta_i$ for all $i=1,\dots,s$. We denote shapes deformed in a  direction $\vec{W} \in C_0^k(D, \R^2)$ by 
\begin{align}
	F_t^{\restr{\vec{W}}{\Delta_i}}(\hat{u}_i)=\{\vec{y} \in D \mid \vec{y}=\vec{x} + t \restr{\vec{W}}{\Delta_i}(\vec{x}) \, \forall\vec{x} \in \hat{u}_i \}
	\label{eq:PerturbationOfIdentity}
\end{align}
 and the vector of shapes (deformed in the $i$-th entry) by $$\hat{\vu}^t_i= (\hat{u}_1,\dots, ,F_t^{\restr{\vec{W}}{\Delta_i}}(\hat{u}_i),\dots,\hat{u}_s).$$
 
\begin{figure}[tb]
	\centering
	\setlength\svgwidth{.55\textwidth}
	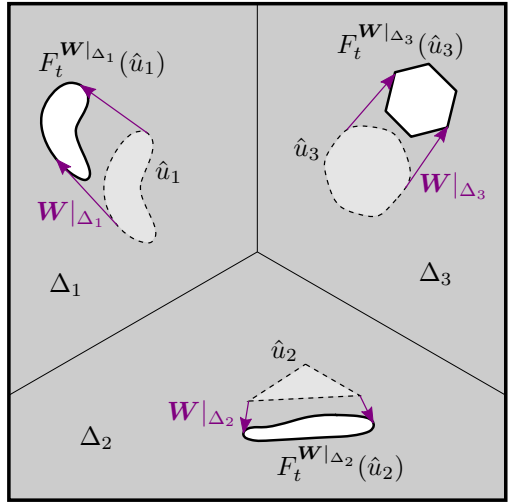
	\caption{Sketch of an admissible partition $\Delta_1,\Delta_2,\Delta_3$ of a domain $D$ containing three shapes $u_1,u_2,u_3$. Each admissible domain $\Delta_i$ is deformed by the vector field $\restr{\vec{W}}{\Delta_i}$.
	}%
	\label{fig:3shapes_drawing_partition}
\end{figure}
\noindent A example of an admissible partition of $D$ and the deformation of the shapes $\hat{\vec{u}}$ can be found in Figure~\ref{fig:3shapes_drawing_partition}.

For $i=1,\dots,s$, the $i$-th partial Eulerian derivative of a function at $\hat{u}_i$ in the direction $\restr{\vec{W}}{\Delta_i}$ is defined by
\begin{equation}
	\label{eulerian_multi_partial}
	\partDeriv \hat{H}(\hat{\vu})[\restr{\vec{W}}{\Delta_i}]\coloneqq  \lim\limits_{t\to 0^+}\frac{\hat{H}(\hat{\vu}^t_i)-\hat{H}(\hat{\vu})}{t}.
\end{equation}
If for all directions $\vec{W} \in C_0^k(D,\R^2)$  and for all $i=1,\dots,s$ the $i$-th partial Eulerian derivative \eqref{eulerian_multi_partial} exists and the mapping 
\begin{equation*}
 C_0^k(D,\R^2) \to \mathbb{R}, \ \vec{W}\mapsto \partDeriv \hat{H}(\hat{\vu})[\restr{\vec{W}}{\Delta_i}]
\end{equation*}
is linear and continuous, the expression $\partDeriv \hat{H}(\hat{\vu})[\restr{\vec{W}}{\Delta_i}]$ is called the $i$-th partial shape derivative of $\hat{H}$ at $\hat{\vu}$ in direction $ \restr{\vec{W}}{\Delta_i}$.
If all these partial shape derivatives exist, then
\begin{equation}
	\label{eulerian_multi_stoch_h}
	\mathrm{d} \hat{H}(\hat{\vu})[\vec{W}]\coloneqq \sum_{i=1}^{s} \partDeriv \hat{H}(\hat{\vu})[\restr{\vec{W}}{\Delta_i}]
\end{equation}
defines the multi-shape derivative of $\hat{H}$ at $\hat{\vu}$ in the direction $\vec{W} \in C_0^k(D,\R^2)$.

We now define the stochastic multi-shape derivative. Let $\hat{J} \colon \hat{M}_s^c \times \Xi \rightarrow \R$ be an objective function that is parametrized with respect to some set $\Xi \subset \R^m$. Let $(\Omega, \mathcal{F}, \pP)$ be a complete probability space, where $\Omega$ is the sample space, $\mathcal{F} \subset 2^{\Omega}$ is the $\sigma$-algebra of events, and $\pP\colon \Omega \rightarrow [0,1]$ is a probability measure. Suppose $\vec{z}\coloneqq\vec{\xi}(\omega)$ is a fixed realization of a random vector $\vec{\xi}\colon \Omega \rightarrow \Xi.$ For $i=1,\dots,s$, the $i$-th partial Eulerian derivative of a the parametrized function $\hat{J}$ at $\hat{u}_i$ (for a fixed realization $\vec{z}$) in the direction $\restr{\vec{W}}{\Delta_i}$ is defined by
\begin{equation}
	\label{eulerian_multi_partial_stoch}
	\partDeriv \hat{J}(\hat{\vu}, \vec{z})[\restr{\vec{W}}{\Delta_i}]\coloneqq  \lim\limits_{t\to 0^+}\frac{\hat{J}(\hat{\vu}^t_i,\vec{z})-\hat{J}(\hat{\vu},\vec{z})}{t}.
\end{equation}
If for all directions $\restr{\vec{W}}{\Delta_i} \in C_0^k(D,\R^2)$ and for all $i=1,\dots,s$  the $i$-th partial Eulerian derivative \eqref{eulerian_multi_partial_stoch} exists and the mapping 
\begin{equation*}
 C_0^k(D,\R^2) \to \mathbb{R}, \ \vec{W}\mapsto \partDeriv \hat{J}(\hat{\vu},\vec{z})[\restr{\vec{W}}{\Delta_i}]
\end{equation*}
is linear and continuous, the expression $\partDeriv \hat{J}(\hat{\vu},\vec{z})[\restr{\vec{W}}{\Delta_i}]$ is called the $i$-th partial shape derivative of $\hat{J}$ at $\hat{\vu}$ (for a fixed realization $\vec{z}$) in direction $\restr{\vec{W}}{\Delta_i}$.
If all these partial shape derivatives exist, then
\begin{equation}
	\label{eulerian_multi_stoch_J}
	\mathrm{d} \hat{J}(\hat{\vu},\vec{z})[\vec{W}]\coloneqq \sum_{i=1}^{s} \partDeriv \hat{J}(\hat{\vu},\vec{z})[\restr{\vec{W}}{\Delta_i}]
\end{equation}
defines the multi-shape derivative of $\hat{J}$ at $\hat{\vu}$  (for a fixed realization $\vec{z}$) in the direction $\vec{W} \in C_0^k(D,\R^2)$. Finally, a stochastic multi-shape derivative is defined as the superposition of 	\eqref{eulerian_multi_stoch_J} with the random vector $\vec{\xi}$, i.e., $\mathrm{d} \hat{J}(\hat{\vu},\vec{\xi})[\vec{W}]$, where for almost every $\omega$ we have
\[
\mathrm{d} \hat{J}(\hat{\vu},\vec{\xi}(\omega))[\vec{W}]\coloneqq \sum_{i=1}^{s} \partDeriv \hat{J}(\hat{\vu},\vec{\xi}(\omega))[\restr{\vec{W}}{\Delta_i}].
\]

%

\section{Application to Navier--Stokes flow} \label{sec:specificProblem}
In this section, we focus on a fluid-mechanical application where viscous energy dissipation is to be minimized with respect to shapes in a connected subset of $\R^2$.   
In Section~\ref{sec:modelFormulation} we introduce the model and the corresponding adjoint equation. Moreover, we provide the detailed shape derivative calculation to the model problem in Section~\ref{subsec:calculationShapeDerivative}.

\subsection{Model formulation} \label{sec:modelFormulation}
Fluid mechanics play a significant role in science and engineering to model the kinematic and dynamic behavior of fluids. 
One part of fluid mechanics is concerned with the evolution of a fluid subject to mass and momentum conservation~\cite{Herwig2006,Danaila2007}. This part can be mathematically described by the Navier--Stokes equations. These equations can be used to model fluid flow of two-dimensional  \cite{Chorin1968,Chen1999,Casas2019,Kuehl2021} to time-dependent three-dimensional problems~\cite{Baensch1991,Alfonsi2014,Bletsos2021}, potentially in real-time with uses in computer graphics, as shown in~\cite{Stam2003,Harris2005,Bridson2008}. 
In this section, we consider steady-state incompressible Navier--Stokes flow with random inputs.
Uncertainty quantification in combination with fluid mechanics has been studied in, e.g., 
\cite{Beck2020a,Benner2020},  
and in the context of shape optimization \cite{Schillings2011a}.

Let $\vu \coloneqq \left( u_1, \ldots, u_s \right) \in M_s^c$ be a vector of shapes belonging to the product shape manifold of piecewise smooth shapes and recall the notation $\hat{\vu} = \left( \hat{u}_1, \ldots, \hat{u}_s \right) \in \hat{M}_s^c$ introduced in the previous section to denote its image in $\R^2$. Consider a hold-all domain $D$ and a non-empty, bounded and connected set $D_{\hat{\vu}} \subset D$ defined such that $\partial D = \Gamma$, $\partial D_{\hat{\vu}} = \Gamma \sqcup \hat{u}_1 \sqcup \cdots \sqcup \hat{u}_s$ and $\hat{\vu} \subset D_{\hat{\vu}}$, where $\Gamma$ is the outer boundary that is fixed and split into two disjoint parts
$	\GD$ and $\GN $ representing a Dirichlet and Neumann boundary, respectively.  
We denote the obstacles outside $D_{\hat{\vu}}$ by $D_{\hat{u}_i}$ such that $\partial D_{\hat{u}_i} = \hat{u}_i$ for all $i$. To avoid trivial solutions in the following problem, additional geometrical constraints have to be added, as discussed in, e.g.,~\cite{Schulz2016}. Here, we implement inequality constraints in our problem description to increase the space of admissible shapes.
For each shape $\hat{u}_i$, we introduce one inequality constraint for its volume, see \eqref{eqn:ModelStokesEquationsVolumeConstraint}, as well as lower and upper bounds for its barycenter, see \eqref{eqn:ModelStokesEquationsBarycenterConstraint}. 
The volume of the domain $D_{\hat{u}_i}$ and the $k$th coordinate of the corresponding barycenter are given by
\begin{subequations}
\label{eq:VolumeBarycenterFormula}
\begin{align}
\vol(\hat{u}_i) &= \intDi{1} = 
 \frac{1}{2} \intui{ \vx \cdot \vn}, \\
   (\bary(\hat{u}_i))_k &= \int_{D_{\hat{u}_i}}\frac{x_k}{\vol(\hat{u}_i)} \d \vx  = \frac{1}{2 \vol(\hat{u}_i)} \intui{ 
x_k^2 n_k}, \label{eq:BarycenterFormula}
\end{align}
\end{subequations}
respectively. Here, the $\vn$ describes the unit outward normal of $D_{\hat{u}_i}$, i.e., the unit inward normal of $D_{\hat{\vu}}$, and $\cdot$ denotes the standard scalar product of two vectors, i.e., $\vec{x} \cdot \vec{n} = \sum_{j=1}^2 x_j n_j$. Further, we require $\vol(\hat{u}_i)>0$ $\forall i=1,\ldots,s$.

The Frobenius inner product for two matrices $\vec{A},\vec{B} \in \R^{2\times 2}$ is denoted by $\vec{A} : \vec{B} = \sum_{j=1}^{2} \sum_{k=1}^{2} A_{j k} B_{j k}$. For a vector field $\vec{v}  \colon \R^2 \rightarrow \R^2$, we define the gradient by $\nabla \vec{v} = (\nabla v_1, \nabla v_2)^\top$, where $\nabla v_j=(\frac{\partial v_j}{\partial x_1},\frac{\partial v_j}{\partial x_2})$, and the Laplacian by $\Delta \vec{v} = (\Delta v_1, \Delta v_2)^\top$, where $\Delta v_j = \sum_{k=1}^2 \frac{\partial^2 v_j}{\partial x_k^2}$. The divergence is defined by $\Divv \vec{v}=\sum_{j=1}^2 \frac{\partial v_j}{\partial x_j}$. The constant $\nu$ represents the fluid viscosity; for simplicity, we do not model this as a random quantity. Additionally, $\vec{f}\colon D \rightarrow \R^2$ is a (deterministic) external force acting on the domain $D$. Now, suppose we have a probability space $(\Omega, \mathcal{F}, \pP)$ and a random vector $\vec{\xi}$ that is compactly supported in $\Xi\subseteq\mathbb{R}^m$. 
Let $\vec{g}\colon D \times \Xi  \rightarrow \R^2$ represent the velocity distribution on $\GD \cup \hat{\vu}$, where for every $\vec{z} \in \Xi$, we have $\vec{g}(\cdot,\vec{z}) \in H^{1/2}(\Gamma_D \cup \hat{\vec{u}})$.
For the random input $\vec{g}(\cdot,\vec{\xi}(\omega))$, we will use the abbreviation $\vec{g}_{\vec{z}}\coloneqq\vec{g}(\cdot,\vec{z})$ for a realization $\vec{z}=\vec{\xi}(\omega) \in \Xi$.

We consider a common shape optimization problem (cf.~\cite{Pironneau1974,Mohammadi2009,Onyshkevych2021})  involving the minimization of dissipated energy, where the physical system is described by the steady-state incompressible Navier--Stokes equations. We modify this problem to include uncertainty in the form of a boundary term $\vec{g}$ as described above; the dissipated energy is to be minimized in expectation. More precisely, we solve the problem
\begin{align}
	\label{eqn:ModelMinimizationProblem}
	\min_{ \hat{\vu}\in \hat{M}_s^c} \quad \left\lbrace \hat{j}(\hat{\vu}) = \frac{\nu}{2} \int_\Omega { \intD{ \nabla \velocity_{\vec{\xi}(\omega)}(\vx) : \nabla \velocity_{\vec{\xi}(\omega)}(\vx) } \d \pP(\omega)} \right\rbrace,
\end{align}
subject to, almost surely,
\begin{subequations}	
\label{eq:ModelStokesEquations}
\begin{align}
	\label{eqn:ModelStokesEquationsMomentum}
	- \nu \Delta \velocity_{\vec{\xi}}(\vx) +  ( \velocity_{\vec{\xi}}(\vx) \cdot \nabla) \velocity_{\vec{\xi}}(\vx) + \nabla \pressure_{\vec{\xi}}(\vx) &= \vec{f}(\vx)  & x \in D_{\hat{\vu}}, \\
	\label{eqn:ModelStokesEquationsContinuity}
	\Div{\velocity_{\vec{\xi}}(\vx) }&= 0 & x \in D_{\hat{\vu}},\\
	\label{eqn:ModelStokesEquationsDirichletBC}
	\velocity_{\vec{\xi}}(\vx) &=\vec{g}_{\vec{\xi}}(\vx)& x \in \GD\cup \hat{\vu}, \\
	\label{eqn:ModelStokesEquationsNeumannBC}
	-\nu \nabla \velocity_{\vec{\xi}}(\vx) \vec{n}(\vx) + \pressure_{\vec{\xi}}(\vx) \vec{n}(\vx) &= \vec{0} & x\in \GN, \\
	\label{eqn:ModelStokesEquationsVolumeConstraint}
	\vol(\hat{u}_i)&\geq \underline{\mathcal{V}}_i  & \forall i=1,\ldots,s ,\\
	\label{eqn:ModelStokesEquationsBarycenterConstraint}
  \underline{\vec{\mathcal{B}}}_{i} \leq \bary(\hat{u}_i) &\leq \overline{\vec{\mathcal{B}}}_{i} & \forall i=1,\ldots,s.
\end{align}
\end{subequations}
Here, $\vec{v}_{\vec{\xi}}$ is the random fluid velocity and $p_{\vec{\xi}}$ is the random pressure such that for almost every $\omega \in \Omega$, we have $\velocity_{\vec{\xi}(\omega)}\colon D_{\hat{\vu}} \rightarrow \R^2$ and $\pressure_{\vec{\xi}(\omega)}\colon D_{\hat{\vu}} \rightarrow \R$.
Together, \eqref{eqn:ModelStokesEquationsMomentum}--\eqref{eqn:ModelStokesEquationsContinuity} are the steady-state incompressible Navier--Stokes equations complemented by the boundary conditions \eqref{eqn:ModelStokesEquationsDirichletBC}--	\eqref{eqn:ModelStokesEquationsNeumannBC}
 with uncertainty induced by the input $\vec{g}_{\vec{\xi}}$. 

Let $V(D_{\hat{\vu}}) = \left\{ \velocity \in H^1( D_{\hat{\vu}}, \R^2)\colon \velocity\vert_{\GD \cup \hat{\vu}} = \vec{0} \right\}$ denote the function space associated to the velocity for a fixed $D_{\hat{\vu}}$.
We assume $\vec{f} \in H^1( D_{\hat{\vu}}, \R^2)$ in view of the higher regularity required for the shape derivative below.
The  weak formulation (see \cite[Chapter 8]{Elman2014}) for a fixed $\vec{z} = \vec{\xi}(\omega)$ is given by: find $\velocity_{\vec{z}} \in H^1( D_{\hat{\vu}}, \R^2)$ and $\pressure_{\vec{z}}   \in L^2( D_{\hat{\vu}})$ such that $\velocity_{\vec{z}} -\vec{g}_{\vec{z}}  \in V(D_{\hat{\vu}})$ and
\begin{subequations}
\label{eq:weak-formulation}
\begin{align}
 \intD{ \nu \nabla \velocity_{\vec{z}}  : \nabla \adjvelocity   +  (( \velocity_{\vec{z}} \cdot \nabla ) \velocity_{\vec{z}}) \cdot \adjvelocity - \pressure_{\vec{z}}  \Div{\adjvelocity } -\vec{f}\cdot \adjvelocity  } &= 0 \quad \forall \adjvelocity  \in V(D_{
 \hat{\vu}}), \\
 \intD{ \adjpressure  \Div{\velocity_{\vec{z}} } } &= 0 \quad \forall \adjpressure \in L^2( D_{
 \hat{\vu}}).
\end{align}
\end{subequations}

\noindent We define $\hat{\vec{h}}\colon \hat{M}_s^c \rightarrow \R^{5s}$ and $\hat{\vec{h}}_i\colon \hat{M}_1^c \rightarrow \R^{5}$ by
\begin{equation*}
\hat{\vec{h}}(\hat{\vu}) =
\begin{pmatrix}
\left[ \mathcal{V}_i - \vol(\hat{u}_i) \right]_{i\in\{1,\dots,s\}} \\[3pt]
\left[ \underline{\vec{\mathcal{B}}}_i -\bary(\hat{u}_i) \right]_{i\in\{1,\dots,s\}} \\[3pt]
\left[ \bary(\hat{u}_i) - \overline{\vec{\mathcal{B}}}_i \right]_{i\in\{1,\dots,s\}}
\end{pmatrix} \quad \text{and} \quad \hat{\vec{h}}_i(\hat{u}_i) = 
\begin{pmatrix}
	\mathcal{V}_i - \vol(\hat{u}_i) 
	\\[3pt]
	\underline{\vec{\mathcal{B}}}_i -\bary(\hat{u}_i) 
	\\[3pt]
	\bary(\hat{u}_i) - \overline{\vec{\mathcal{B}}}_i 
\end{pmatrix}
\end{equation*}
and the parametrized objective function by
$\hat{J}(\hat{\vu},\vec{z})\coloneqq\frac{\nu}{2}  \int_{D_{\hat{\vu}}} \nabla \velocity_{\vec{z}}(\vx): \nabla \velocity_{\vec{z}}(\vx) \d \vx$. As we will rely on the stochastic augmented Lagrangian method in computations, we will also need the parametrized augmented Lagrangian, defined for $\mu >0$ by 
\begin{align}
\label{eqn:augLagrangeFunction}
\begin{aligned}
	\hat{L}_{A}(\hat{\vu}, \vec{\lambda}, \vec{z};\mu) 
	=\;&   \hat{J}(\hat{\vu},\vec{z}) + \intD{ \nu \nabla \velocity_{\vec{z}} : \nabla \adjvelocity_{\vec{z}}  +(( \velocity_{\vec{z}}  \cdot \nabla) \velocity_{\vec{z}} ) \cdot \adjvelocity_{\vec{z}}}\\
	&  + \intD{ - \pressure_{\vec{z}}  \Div{\adjvelocity_{\vec{z}} } -\vec{f} \cdot \adjvelocity_{\vec{z}} + \adjpressure_{\vec{z}}  \Div{\velocity_{\vec{z}} } } \\
	& + \frac{\mu}{2} \left\lVert \max\left(0, \hat{\vec{h}}(\hat{\vu}) + \frac{\vec{\lambda}}{\mu} \right)\right\rVert_2^2 - \frac{\|\vec{\lambda}\|_2^2}{2 \mu},
\end{aligned}
\end{align}
where  $\velocity_{\vec{z}} \in H^1( D_{\hat{\vu}}, \R^2)$ and $\pressure_{\vec{z}}   \in L^2( D_{\hat{\vu}})$ solve \eqref{eq:weak-formulation} and $\adjvelocity_{\vec{z}} \in V(D_{\hat{\vu}})$ and $\adjpressure_{\vec{z}} \in L^2(D_{\hat{\vu}})$ solve the (weak form of the) adjoint equation (cf., e.g.,~\cite[p.~114]{hinze2002}) for a fixed $\vec{z} \in \Xi$:
\begin{subequations}
\label{eqn:ModelStokesEquationsMomentumAdj}	
\begin{align}
	\begin{aligned}
\intD{ \nu \nabla \tilde{\adjvelocity} : (\nabla \velocity_{\vec{z}}+\nabla\adjvelocity_{\vec{z}}) +((\tilde{\adjvelocity}\cdot\nabla )\velocity_{\vec{z}}) \cdot \adjvelocity_{\vec{z}} \quad & \\
+((\velocity_{\vec{z}}\cdot\nabla )\tilde{\adjvelocity}) \cdot \adjvelocity_{\vec{z}} + \adjpressure_{\vec{z}} \Div{\tilde{\vec{\varphi}}} &} = 0 \quad \forall \tilde{\adjvelocity} \in V(D_{\hat{\vu}}),
	\end{aligned} \\ 
\intD{ -\Div{\adjvelocity_{\vec{z}}} \,\tilde{\adjpressure}} = 0 \quad \forall \tilde{\adjpressure} \in L^2(D_{\hat{\vu}}).
\end{align}
\end{subequations}

\subsection{Calculation of the stochastic multi-shape derivative}\label{subsec:calculationShapeDerivative}
To compute the multi-shape derivative for problem~\eqref{eqn:ModelMinimizationProblem}--\eqref{eq:ModelStokesEquations}, let
$\hat{L}_A\colon \hat{M}^c_s \times \R^n \times \Xi \rightarrow \R$ and $\hat{\vec{h}} \colon \hat{M}^c_s \rightarrow \R$ be the functions defined by $\hat{L}_A(\hat{\vu}, \vec{\lambda}, \vec{z}; \mu)={L}_A(\mathbb{U}({\vu}), \vec{\lambda}, \vec{z}; \mu)$ and $\hat{\vec{h}}(\hat{\vec{\vu}}) = \vec{h}(\mathbb{U}(\vu))$ for all $\vu \in M_c^s$. We consider vector fields from the set $\mathcal{W}(D_{\hat{\vu}})\coloneqq \{ \vec{W} \in H^1(D_{\hat{\vu}}, \R^2)\colon \vec{W}|_{\Gamma}=0 \}$; the condition $\vec{W}|_{\Gamma}=0$ ensures that the outer boundary of the hold-all domain does not move.

	There are a lot of options to calculate the shape derivative of shape functionals. An overview about the min-max approach \cite{Delfour-Zolesio-2001}, the chain rule approach \cite{SokoZol}, the Lagrange method of C\'{e}a \cite{Cea-RAIRO} and the rearrangement method \cite{Ito-Kunisch-Peichl} is given in \cite{Sturm}. 
	In addition, there are so-called material-free approaches available (see, e.g.,  \cite{Sturm2013}).
	In this paper, we focus on a material derivative approach. 
	
To calculate the shape derivative of \eqref{eqn:augLagrangeFunction}, we need to define the multi-material derivative and to derive several expressions, e.g., the material derivative of vector-valued functions.
Then, we will be able to calculate the shape derivative of the augmented Lagrangian.

\subsubsection{Necessary formulas regarding the shape derivative}

\paragraph{Multi-material derivative}
A definition of the material derivative in the setting of product spaces needs to be specified; we will call it the multi-material derivative in the following.
For the moment, we consider the setting from Section~\ref{subsec:multi-shape-derivative},  i.e., that $D$ is a domain that contains $\hat{\vec{u}}\in \hat{M}_s^c$ and is partitioned into subsets $\Delta_1,\dots,\Delta_s$ such that $\hat{u}_i\subset \Delta_i$ for all $i=1,\dots,s$. 
For  $\vec{W} \in C_0^k(D, \R^2)$, we define 
\begin{align*}
\Delta_i^t\coloneqq F_t^{\restr{\vec{W}}{\Delta_i}}(\Delta_i)=\{\vec{y} \in D \mid \vec{y}=\vec{x} + t \restr{\vec{W}}{\Delta_i}(\vec{x}) \, \forall\vec{x} \in \Delta_i \}
\end{align*}
for all $i=1,\dots, s$ and $t\geq 0$.
The domain $D$ can then be seen as dependent on the transformations $F_t$. Therefore, it is convenient to define
$D^t\coloneqq \left(\cup_{i=1}^s \Delta_i^t \right) \cup  \partial \Delta^t$ with $\partial \Delta^t\coloneqq \cup_{i=1}^s \partial\Delta_i^t\setminus \partial D$ and $\partial\Delta\coloneqq\partial\Delta^0$. Thanks to this partition, we notice that a point $\vec{x}\in D$ 
can be related to the corresponding $\vec{x}\in \Delta_i$ or $\vec{x}\in \partial \Delta$ 
for all $i=1,\dots,s$. 
Now, we can formulate the definition of the multi-material  derivative.
\begin{definition}
	\label{def:multiMaterialDerivative}
	Let $p\colon D \to \mathbb{R}$. 
	For  $\vec{W} \in C_0^k(D, \R^2)$, we consider $D^t\coloneqq \left(\cup_{i=1}^s \Delta_i^t \right) \cup  \partial \Delta^t$.
	Moreover, let $\{p^t \colon D^t\to\R,\, t\leq T\}$, $\{p_b^t\coloneqq\restr{p^t}{\partial\Delta^t}  \colon \partial\Delta^t\to\R,\, t\leq T\}$, and $\{p_i^t\coloneqq\restr{p^t}{\Delta_i^t}  \colon \Delta_i^t\to\R,\, t\leq T\}$ with $T>0$ denote families of mappings for all $i=1,\dots, s$. 
	We define 
		\begin{equation*}
			\partMatDeriv p(\vec{x})  \coloneqq
		\begin{cases}
		\lim\limits_{t\to 0^+}\frac{\left(p_i^t\circ F_t^{\restr{\vec{W}}{\Delta_i}}\right)(\vec{x})-p_i^0(\vec{x})}{t} =\frac{\d^+}{\d t}\left(p_i^t\circ F_t^{\restr{\vec{W}}{\Delta_i}}\right)(\vec{x})\,\rule[-2.5mm]{.1mm}{6mm}_{\hspace{1mm}t=0} & \text{ if }\vec{x}\in \Delta_i,\\[8pt]
			0  &\text{ else}
		\end{cases}
		\end{equation*}
		for all $i=1,\dots, s$ and 
			\begin{equation*}
				\d_{m_{s+1}}p(\vec{x})  \coloneqq
			\begin{cases}
	 \lim\limits_{t\to 0^+}\frac{\left(p_b^t\circ F_t^{\restr{\vec{W}}{\partial\Delta}}\right)(\vec{x})-p_b^0(\vec{x})}{t} =\frac{\d^+}{\d t}\left(p_b^t\circ F_t^{\restr{\vec{W}}{\partial\Delta}}\right)(\vec{x})\,\rule[-2.5mm]{.1mm}{6mm}_{\hspace{1mm}t=0} & \text{ if }\vec{x}\in\partial \Delta,\\[8pt]
	 0 &\text{ else.}
			\end{cases}
	\end{equation*}
	We call $\partMatDeriv p(\vec{x})$ the $i$-th partial material derivative of $p$ at $\vec{x}\in D$.
	The multi-material derivative of  $p$ at $\vec{x}\in D$ is denoted by $\matDeriv p(\vec{x})$ or $\dot{p}(\vec{x})$ and given by 
	\begin{equation*}
	\matDeriv p(\vec{x})= 
	\sum_{i=1}^{s+1} \partMatDeriv p(\vec{x}) .
	\end{equation*}
\end{definition}

\paragraph{Material derivative of vector-valued functions}
For the moment, we consider the general setting in $\R^d$ with $\tilde{D} \subset \R^d$. The material derivative of a sufficiently smooth function $\vec{p}=(p_1 , \ldots, p_d)^\top \colon \tilde{D
} \times \R \rightarrow \R^d$ is denoted by $\matDeriv(\vec{p})$ or $\dot{\vec{p}}$, where the second argument denotes the size of perturbation $t$ as in \eqref{eq:PerturbationOfIdentity}. Similar to \cite{Berggren2009}, we suppress the second argument in case there is no possibility of confusion. We will also need a generalization of~\cite[equation~(25)]{Berggren2009} for vector-valued functions: the material derivative of $\nabla \vec{p}$ in the direction of a sufficiently smooth vector field $\vec{W}\in C_0^k(\tilde{D},\R^d)$ is given by
\begin{align}
	\matDeriv(\nabla \vec{p}) 
=  \nabla \matDeriv({\vec{p}}) - \nabla \vec{p} \nabla \vec{W}.
	\label{eqn:MaterialDerivativeGradientVectorValued}
\end{align}
Given a second function $\vec{q}=(q_1 , \ldots, q_d)^\top\colon \tilde{D} \times \R \rightarrow \R^d$, we have the following identity 
\begin{align}
	\nabla \vec{p} \nabla \vec{q} : \mathbb{I} = \sum\limits_{k=1}^d\sum\limits_{j=1}^d\sum\limits_{i=1}^d \frac{\partial p_i}{\partial x_j} \frac{\partial q_j}{\partial x_k} \delta_{i k}= \sum\limits_{k=1}^d\sum\limits_{j=1}^d \frac{\partial p_k}{\partial x_j} \frac{\partial q_j}{\partial x_k} = {\nabla \vec{p}}^\top : \nabla \vec{q},
	\label{eqn:MaterialDerivativeDivergence}
\end{align}
where $\delta_{i k}$ is the Kronecker delta and $\mathbb{I}\in \mathbb{R}^{d\times d}$ the identity matrix. 
		
\paragraph{Quotient shape derivative rule}
We consider again the general setting in $\R^d$ with $\tilde{D} \subset \R^d$.
The shape derivative of the barycenter constraint term \eqref{eq:BarycenterFormula} is required. Let $q\colon \tilde{D} \times \R \rightarrow \R$  and $\vec{W}\in C_0^k(\tilde{D},\R^d)$ be sufficiently smooth. 
We consider the domain integral $Q(t)\coloneqq\int_{\tilde{D}} q(\vec{x},t) \dx$ and assume $Q(t) \neq 0$ for all $t\geq 0$ in a neighborhood of zero.
 Adapting the standard definition of the shape derivative of a domain integral  $Q(t)$ (cf.~\cite{Haslinger2003,Berggren2009}) yields
\begin{align}
	\begin{aligned}
	\label{eqn:ShapeDerivativeQuotientRule}
	\mathrm{d} Q^{-1}(0)[\vec{W}] 
	&
	= \frac{\mathrm{d}}{\mathrm{d}t} \left. \frac{1}{\int_{\tilde{D}} q(F_t^{\vec{W}}(\vec{x}),t) \,\dx} \right|_{t=0^+} 
	= - \frac{\frac{\mathrm{d}}{\mathrm{d}t} \left. \int_{\tilde{D}} q(F_t^{\vec{W}}(\vec{x}),t) \,\dx \right|_{t=0^+}}{\left( \int_{\tilde{D}} q(\vec{x},0) \,\dx \right)^2} \\
	&= - \frac{\mathrm{d} Q(0)[\vec{W}] }{Q(\vec{x},0)^2}.
	\end{aligned}
\end{align}

\subsubsection{Shape derivative of the augmented Lagrangian} 
We return to our example and now split the parametrized augmented Lagrangian $\hat{L}_A\colon \hat{M}_c^s \times \R^{5s} \times \Xi \rightarrow \R$ 
into separate terms via
\begin{align*}
	\hat{L}_{A}(\hat{\vu}, \vec{\lambda}, \vec{z};\mu) = \sum_{\chi=1}^{8} \hat{L}_{\chi}(\hat{\vu}, \vec{\lambda}, \vec{z};\mu),
\end{align*}
where
\begin{align*}
	\hat{L}_{1}(\hat{\vu}, \vec{\lambda}, \vec{z};\mu) &= \intD{\frac{\nu}{2} \, \nabla \velocity : \nabla \velocity}, &
	\hat{L}_{2}(\hat{\vu}, \vec{\lambda}, \vec{z};\mu) &= \intD{\nu \nabla \velocity : \nabla \adjvelocity}, \\
	\hat{L}_{3}(\hat{\vu}, \vec{\lambda}, \vec{z};\mu) &= \intD{(( \velocity_{\vec{z}}  \cdot \nabla) \velocity_{\vec{z}} ) \cdot \adjvelocity_{\vec{z}}}, &
	\hat{L}_{4}(\hat{\vu}, \vec{\lambda}, \vec{z};\mu) &= \intD{ - \pressure \Div{\adjvelocity}}, \\
	\hat{L}_{5}(\hat{\vu}, \vec{\lambda}, \vec{z};\mu) &= \intD{-\vec{f}\cdot \adjvelocity}, &
	\hat{L}_{6}(\hat{\vu}, \vec{\lambda}, \vec{z};\mu) &= \intD{\adjpressure \Div{\velocity}}, \\
	\hat{L}_{7}(\hat{\vu}, \vec{\lambda}, \vec{z};\mu) &= \frac{\mu}{2}\left \lVert \max \left(0,\hat{\vec{h}}(\hat{\vec{u}}) + \frac{\vec{\lambda}}{\mu} \right)\right\rVert_2^2, \quad \text{and} &
	\hat{L}_{8}(\hat{\vu}, \vec{\lambda}, \vec{z};\mu) &= - \frac{\|\vec{\lambda}\|_2^2}{2 \mu}.
\end{align*}

To compute the shape derivative of each $\hat{L}_{\chi}$ in direction \mbox{$\vec{W}\in \mathcal{W}(D_{\hat{\vu}})$} and for a fixed realization, we use the definition of the material derivative and shape derivative as described in \cite[Definition 1 and 2]{Berggren2009} with respect to $\vec{W}$, the concept of partial shape derivatives and partial material derivatives described in Section~\ref{subsec:multi-shape-derivative} and Definition~\ref{def:multiMaterialDerivative}, and the linearity of the shape derivative. This gives 
\begin{equation}
	\label{eqn:ShapeDerivativeSplitting}
	\begin{split}
&	\Deriv \hat{L}_{A}(\hat{\vu}, \vec{\lambda},\vec{z};\mu)[\vec{W}] \\&= \sum_{i=1}^{s} \partDeriv \hat{L}_A(\hat{\vu}, \vec{\lambda},\vec{z};\mu)[\restr{\vec{W}}{\Delta_i}] = \sum_{i=1}^{s} \sum_{\chi=1}^{8} \partDeriv \hat{L}_{\chi}(\hat{\vu}, \vec{\lambda},\vec{z};\mu)[\restr{\vec{W}}{\Delta_i}].
	\end{split}
\end{equation}	

Using the rules of the material derivative in \cite[section 5]{Berggren2009} and \eqref{eqn:MaterialDerivativeGradientVectorValued}, as well as $\matDeriv(\nu)=0$, we get 
\begin{subequations}
	\begin{align}
		\label{eqn:shapeDerivL1}
		\begin{aligned}
			&\partDeriv \hat{L}_1(\hat{\vu}, \vec{\lambda},\vec{z};\mu)\left[\restr{\vec{W}}{\Delta_i}\right]=\partDeriv \left( \intDeltai{\frac{\nu}{2} \, \nabla \velocity_{\vec{z}} : \nabla \velocity_{\vec{z}}} \right) \left[\restr{\vec{W}}{\Delta_i}\right] \\
			&= \intDeltai{ \frac{\nu}{2} \left( \partMatDeriv ( \nabla \velocity_{\vec{z}} ) : \nabla \velocity_{\vec{z}} + \nabla \velocity_{\vec{z}} : \partMatDeriv (\nabla \velocity_{\vec{z}}) \right) 
				\\
				&	\phantom{= \int_{Delta_i} }+ \Div{\restr{\vec{W}}{\Delta_i}} \left( \frac{\nu}{2} \, \nabla \velocity_{\vec{z}} : \nabla \velocity_{\vec{z}} \right)  } \\
			&= \intDeltai{ \nu \nabla \partMatDeriv (\velocity_{\vec{z}}) : \nabla \velocity_{\vec{z}} - \nu \left(\nabla \velocity_{\vec{z}} \nabla \restr{\vec{W}}{\Delta_i} \right) : \nabla \velocity_{\vec{z}} \\
			&	\phantom{= \int_{Delta_i} }+ \Div{\restr{\vec{W}}{\Delta_i}} \left( \frac{\nu}{2} \, \nabla \velocity_{\vec{z}} : \nabla \velocity_{\vec{z}} \right)  }.
		\end{aligned}
	\end{align}
	Performing an analogous calculation for $\partDeriv \hat{L}_2(\hat{\vu}, \vec{\lambda},\vec{z};\mu)\left[\restr{\vec{W}}{\Delta_i}\right]$ yields
	\begin{align}
		\label{eqn:shapeDerivL2}
		\begin{aligned}
			&\partDeriv \hat{L}_2(\hat{\vu}, \vec{\lambda},\vec{z};\mu)\left[\restr{\vec{W}}{\Delta_i}\right]=\partDeriv \left( \intDeltai{ \nu \nabla \velocity_{\vec{z}} : \nabla \adjvelocity_{\vec{z}}} \right) \left[\restr{\vec{W}}{\Delta_i}\right] \\
			&= \intDeltai{ \nu \nabla \partMatDeriv (\velocity_{\vec{z}}) : \nabla \adjvelocity_{\vec{z}} - \nu \left( \nabla \velocity_{\vec{z}} \nabla \restr{\vec{W}}{\Delta_i} \right) : \nabla \adjvelocity_{\vec{z}} + \nu \nabla \velocity_{\vec{z}} : \nabla \partMatDeriv (\adjvelocity_{\vec{z}}) \\
				&\hphantom{=\int_{D_{\vu}}}\ - \nu \nabla \velocity_{\vec{z}} : \left( \nabla \adjvelocity_{\vec{z}} \nabla \restr{\vec{W}}{\Delta_i} \right) + \Div{\restr{\vec{W}}{\Delta_i}} \left( \nu \nabla \velocity_{\vec{z}} : \nabla \adjvelocity_{\vec{z}} \right) }.
		\end{aligned}
	\end{align}
	The shape derivative of $\hat{L}_3$ is given by
	\begin{align}
		\label{eqn:shapeDerivL3}
		\begin{aligned}
			&\partDeriv \hat{L}_3(\hat{\vu}, \vec{\lambda},\vec{z};\mu)\left[\restr{\vec{W}}{\Delta_i}\right]=\partDeriv \left( \intDeltai{(( \velocity_{\vec{z}}  \cdot \nabla) \velocity_{\vec{z}} ) \cdot \adjvelocity_{\vec{z}}} \right) \left[\restr{\vec{W}}{\Delta_i}\right] \\
			&= \intDeltai{ \partMatDeriv ( (\velocity_{\vec{z}}  \cdot \nabla) \velocity_{\vec{z}} ) \cdot \adjvelocity_{\vec{z}} + ( (\velocity_{\vec{z}}  \cdot \nabla) \velocity_{\vec{z}} ) \cdot \partMatDeriv (\adjvelocity_{\vec{z}}) \\
			&\hphantom{=\int_{D_{\vu}}}+ \Div{\restr{\vec{W}}{\Delta_i}} \left(  ((\velocity_{\vec{z}}  \cdot \nabla) \velocity_{\vec{z}} ) \cdot \adjvelocity_{\vec{z}} \right) } \\
			&= \intDeltai{  ( (\partMatDeriv (\velocity_{\vec{z}})  \cdot \nabla) \velocity_{\vec{z}} ) \cdot \adjvelocity_{\vec{z}} + ( (\velocity_{\vec{z}}  \cdot \nabla) \partMatDeriv (\velocity_{\vec{z}}) ) \cdot \adjvelocity_{\vec{z}} \\
			&\hphantom{=\int_{D_{\vu}}}	
			-	 ((\nabla \restr{\vec{W}}{\Delta_i} \velocity_{\vec{z}} \cdot \nabla) \velocity_{\vec{z}}) \cdot \adjvelocity_{\vec{z}} 
			+ ( (\velocity_{\vec{z}}  \cdot \nabla) \velocity_{\vec{z}} ) \cdot \partMatDeriv (\adjvelocity_{\vec{z}})  
			\\ &\hphantom{=\int_{D_{\vu}}}\ + \Div{\restr{\vec{W}}{\Delta_i}} \left(  ((\velocity_{\vec{z}}  \cdot \nabla) \velocity_{\vec{z}} ) \cdot \adjvelocity_{\vec{z}} \right) }.
		\end{aligned}
	\end{align}
	Using \eqref{eqn:MaterialDerivativeDivergence} and $\partMatDeriv\mathbb{I}=\vec{0}$, we get the shape derivative 
	\begin{align}
		\label{eqn:shapeDerivL4}
		\begin{aligned}
			&\partDeriv \hat{L}
			_4(\hat{\vu}, \vec{\lambda},\vec{z};\mu)\left[\restr{\vec{W}}{\Delta_i}\right]=\partDeriv \left( \intDeltai{ - \pressure_{\vec{z}} \Div{\adjvelocity_{\vec{z}}} } \right) \left[\restr{\vec{W}}{\Delta_i}\right] \\
			&= \intDeltai{ -\partMatDeriv (\pressure_{\vec{z}}) \Div{\adjvelocity_{\vec{z}}} - \pressure_{\vec{z}} \partMatDeriv (\nabla \adjvelocity_{\vec{z}} : \mathbb{I}) + \Div{\restr{\vec{W}}{\Delta_i}} \left( - \pressure_{\vec{z}} \Div{\adjvelocity_{\vec{z}}} \right) } \\
			&= \intDeltai{ - \partMatDeriv (\pressure_{\vec{z}}) \Div{\adjvelocity_{\vec{z}}} - \pressure_{\vec{z}} \left(\nabla \partMatDeriv (\adjvelocity_{\vec{z}}) - \nabla \adjvelocity_{\vec{z}} \nabla \restr{\vec{W}}{\Delta_i}\right) : \mathbb{I}\\
				&\hphantom{=\int_{D_{\vu}}} + \Div{\restr{\vec{W}}{\Delta_i}} \left( - \pressure_{\vec{z}} \Div{\adjvelocity_{\vec{z}}} \right) } \\
			&= \intDeltai{ - \partMatDeriv (\pressure_{\vec{z}}) \Div{\adjvelocity_{\vec{z}}} - \pressure_{\vec{z}} \Div{\partMatDeriv (\adjvelocity_{\vec{z}})} + \pressure_{\vec{z}} {\nabla \adjvelocity_{\vec{z}}}^\top : \nabla \restr{\vec{W}}{\Delta_i} \\
			&\hphantom{=\int_{D_{\vu}}} + \Div{\restr{\vec{W}}{\Delta_i}} \left( - \pressure_{\vec{z}} \Div{\adjvelocity_{\vec{z}}} \right) }.
		\end{aligned}
	\end{align}
	Using $\dot{\vec{f}}=\nabla \vec{f} \restr{\vec{W}}{\Delta_i}$ since $\vec{f}$ is not dependent on the shape, see e.g.~\cite[Chapter 4, Example 2]{Berggren2009}, the shape derivative $\partDeriv \hat{L}_5(\hat{\vu}, \vec{\lambda},\vec{z};\mu)\left[\restr{\vec{W}}{\Delta_i}\right]$ results in
	\begin{align}
		\label{eqn:shapeDerivL5}
		\begin{aligned}
			&\partDeriv \hat{L}_5(\hat{\vu}, \vec{\lambda},\vec{z};\mu)\left[\restr{\vec{W}}{\Delta_i}\right]=\partDeriv \left( \intDeltai{ - \vec{f} \cdot \adjvelocity_{\vec{z}} } \right) \left[\restr{\vec{W}}{\Delta_i}\right] \\
			&= \intDeltai{ - \partMatDeriv \left( \vec{f} \right) \cdot \adjvelocity_{\vec{z}}  - \vec{f} \cdot \partMatDeriv (\adjvelocity_{\vec{z}} ) + \Div{\restr{\vec{W}}{\Delta_i}}  \left( - \vec{f} \cdot \adjvelocity_{\vec{z}} \right) } \\
			&= \intDeltai{ - (\nabla \vec{f} \restr{\vec{W}}{\Delta_i}) \cdot \adjvelocity_{\vec{z}} - \vec{f} \cdot \partMatDeriv (\adjvelocity_{\vec{z}}) + \Div{\restr{\vec{W}}{\Delta_i}}  \left( - \vec{f} \cdot \adjvelocity_{\vec{z}} \right) }.
		\end{aligned}
	\end{align}
	The shape derivative  of $\hat{L}_6$ can be calculated analogously to the one of $\hat{L}_4$. We get
	\begin{align}
		\label{eqn:shapeDerivL6}
		\begin{aligned}
			&\partDeriv \hat{L}_6(\hat{\vu}, \vec{\lambda},\vec{z};\mu)\left[\restr{\vec{W}}{\Delta_i}\right]=\partDeriv \left( \intDeltai{ \adjpressure_{\vec{z}} \Div{\velocity_{\vec{z}}} } \right) \left[\restr{\vec{W}}{\Delta_i}\right]	\\
			&= \intDeltai{ \partMatDeriv (\adjpressure_{\vec{z}}) \Div{\velocity_{\vec{z}}} + \adjpressure_{\vec{z}} \Div{\partMatDeriv (\velocity_{\vec{z}})} - \adjpressure_{\vec{z}} {\nabla \velocity_{\vec{z}}}^\top : \nabla \restr{\vec{W}}{\Delta_i} \\
				&\quad \qquad
				+ \Div{\restr{\vec{W}}{\Delta_i}} \left( \adjpressure_{\vec{z}} \Div{\velocity_{\vec{z}}} \right) }.
		\end{aligned}
	\end{align}
\end{subequations}


\noindent Choosing a $\widetilde{\vec{W}} \in H^1(D_{\hat{u}_i}, \R^2)$ for the computation of $\partDeriv \hat{L}_7(\hat{\vu}, \vec{\lambda},\vec{z};\mu)\left[\restr{\vec{W}}{\Delta_i}\right]$ yields
\begin{subequations}
	\begin{align*}
		\begin{aligned}
		&	\partDeriv\! \left( \frac{\mu}{2} \left \lVert \max \left(0,\hat{\vec{h}}(\hat{\vec{u}}) + \frac{\vec{\lambda}}{\mu} \right)\right\rVert_2^2\right)\! \left[\widetilde{\vec{W}}\right] \\
			&\quad = \mu\! \left( \left(\hat{\vec{h}}(\hat{\vec{u}}) + \frac{\vec{\lambda}}{\mu} \right) -  \max \left(0,\hat{\vec{h}}(\hat{\vec{u}}) + \frac{\vec{\lambda}}{\mu} \right) \right) \partDeriv\! \left(\hat{\vec{h}}(\hat{\vec{u}}) + \frac{\vec{\lambda}}{\mu} \right)\! \left[\widetilde{\vec{W}}\right] \\
			&\quad = \mu \left( \left(\hat{\vec{h}}(\hat{\vec{u}}) + \frac{\vec{\lambda}}{\mu} \right) -  \max \left(0,\hat{\vec{h}}(\hat{\vec{u}}) + \frac{\vec{\lambda}}{\mu} \right) \right) \partDeriv (\hat{\vec{h}}(\hat{\vu})) \!\left[\widetilde{\vec{W}}\right].
		\end{aligned}
	\end{align*}
	It is required to calculate $\partDeriv \hat{\vec{h}}(\hat{\vu}) \!\left[\widetilde{\vec{W}}\right]$, where we first note that the partial shape derivative of any constraints in $\hat{\vec{h}}(\hat{\vu})$ is zero, except for $\hat{\vec{h}}_i(\hat{u}_i)$, and insert \eqref{eq:VolumeBarycenterFormula}:

	\begin{align*}
		\partDeriv \hat{\vec{h}}_i(\hat{u}_i) \!\left[\widetilde{\vec{W}}\right] = \partDeriv \begin{pmatrix}
			\mathcal{V}_i- \vol(\hat{u}_i) 
			\\[3pt]
			\underline{\vec{\mathcal{B}}}_i -\bary(\hat{u}_i) 
			\\[3pt]
			\bary(\hat{u}_i) - \overline{\vec{\mathcal{B}}}_i 
		\end{pmatrix} \!\left[\widetilde{\vec{W}}\right] &
		= \partDeriv \begin{pmatrix}
			-  \intDi{1} 
			\\[3pt]
			- \frac{\intDi{\vec{x}} }{\intDi{1} } 
			\\[3pt]
			\frac{\intDi{\vec{x}} }{\intDi{1} } 
		\end{pmatrix} \!\left[\widetilde{\vec{W}}\right].
	\end{align*}
	Then, using the divergence theorem with the unit inward normal $\vec{n}$ of $D_{\hat{u}_i}$, which is equivalent to the unit outward normal of $D_{\hat{\vu}}$, gives the identity
	$$	\partDeriv \left( \intDi{1} \right) \!\left[\widetilde{\vec{W}}\right] = 
	\intDi{\Div{\widetilde{\vec{W}}}} = \intui{-\widetilde{\vec{W}} \cdot \vec{n}}.$$
Similar to~\cite[section 4, example 1]{Berggren2009}, using the identity above and
	applying the quotient rule of the shape derivative given in \eqref{eqn:ShapeDerivativeQuotientRule} with $\intDi{1}>0$ $\forall i=1,\ldots,s$, see \eqref{eq:VolumeBarycenterFormula}, we get
	\begin{align*}
		\partDeriv &\left( \frac{\intDi{\vec{x}} }{\intDi{1} } \right) \!\left[\widetilde{\vec{W}}\right] \\
		&= \frac{1}{\intDi{1}} \partDeriv \left( \intDi{\vec{x}} \right) \!\left[\widetilde{\vec{W}}\right] + \partDeriv \left( \frac{1}{\intDi{1}} \right) \!\left[\widetilde{\vec{W}}\right] \intDi{\vec{x}} \\
		&= \frac{1}{\intDi{1}} \intDi{\widetilde{\vec{W}} \cdot \nabla \vec{x} + \vec{x} \Div{\widetilde{\vec{W}}}} - \frac{\intDi{\vec{x}}}{\left(\intDi{1}\right)^2} \intDi{\Div{\widetilde{\vec{W}}}}\\ 
		&= \frac{1}{\vol(\hat{u}_i)} \intui{ \vec{x} \left(\widetilde{\vec{W}} \cdot (-\vec{n}) \right) } - \frac{\bary(\hat{u}_i)}{\vol(\hat{u}_i)} \intui{\widetilde{\vec{W}} \cdot \left(-\vec{n}\right)} \\
		&= \frac{1}{\vol(\hat{u}_i)} \intui{- \left( \vec{x} - \bary(\hat{u}_i) \right) \left( \widetilde{\vec{W}} \cdot \vec{n} \right) }.
	\end{align*}
	The direction $\widetilde{\vec{W}}$ on $u_i$ is equivalent to the direction $\vec{W} \in \mathcal{W}(D_{\hat{\vu}})$ on $u_i$, therefore the shape derivative $\partDeriv \hat{L}_7(\hat{\vu}, \vec{\lambda},\vec{z};\mu)\left[\restr{\vec{W}}{\Delta_i}\right]$ follows as
	\begin{align}
		\label{eqn:shapeDerivL7}
		\begin{aligned}
		\partDeriv \hat{L}_7(\hat{\vu}, \vec{\lambda},\vec{z};\mu)\left[\restr{\vec{W}}{\Delta_i}\right] &=
		\mu \left( \left(\hat{\vec{h}}_i(\hat{u}_i) + \frac{\vec{\lambda}_i}{\mu} \right) - \max \left(0, \hat{\vec{h}}_i(\hat{u}_i) + \frac{\vec{\lambda}_i}{\mu} \right) \right)  \\
		&\qquad \cdot \begin{pmatrix}
				\intui{\restr{\vec{W}}{\Delta_i} \cdot \vec{n}}
				\\[3pt]
				\frac{1}{\vol(\hat{u}_i)} \intui{ \left( \vec{x} - \bary(\hat{u}_i) \right) \left( \restr{\vec{W}}{\Delta_i} \cdot \vec{n} \right) } 
				\\[3pt]
				\frac{1}{\vol(\hat{u}_i)} \intui{ \left( \bary(\hat{u}_i)- \vec{x} \right) \left( \restr{\vec{W}}{\Delta_i} \cdot \vec{n} \right) } 
			\end{pmatrix}
		\end{aligned}
	\end{align}
where $\vec{\lambda}_i$ denotes the Lagrange multipliers related to the contraints $\hat{\vec{h}}_i(\hat{u}_i)$.
	Lastly, the shape derivative $\partDeriv \hat{L}_8(\hat{\vec{\vu}}, \vec{\lambda},\vec{z};\mu)\left[\restr{\vec{W}}{\Delta_i}\right]$ is
	\begin{align}
		\label{eqn:shapeDerivL8}
		\partDeriv \hat{L}_8(\hat{\vec{u}}, \vec{\lambda},\vec{z};\mu)\left[\restr{\vec{W}}{\Delta_i}\right]=\partDeriv \left( \frac{\|\vec{\lambda}\|_2^2}{2 \mu} \right) \left[\restr{\vec{W}}{\Delta_i}\right] = 0.
	\end{align}
\end{subequations}



In order to obtain the shape derivative of $\hat{L}_A$ in direction $\vec{W} \in \mathcal{W}(D_{\hat{\vu}})$ we we now combine the partial shape derivatives $\partDeriv \hat{L}_\chi(\hat{\vu}, \vec{\lambda}, \vec{z};\mu)[\restr{\vec{W}}{\Delta_i}]$, $\chi=1,\ldots,8$ and $i=1,\ldots,s$, according to \eqref{eqn:ShapeDerivativeSplitting}. 
Moreover, we note that $\partial \Delta$ describes a null set as a one-dimensional subspace of $\R^2$ and therefore $\sum_{i=1}^s \intDeltai{\ldots}=\intD{\ldots}$. In summary, we get
\begin{align*}
&\hspace*{-.2cm}	\Deriv  \hat{L}_{A}(\hat{\vu}, \vec{\lambda},\vec{z};\mu)[\vec{W}]\\
	=& \int_{D_{\hat{\vu}}} \nu \nabla \dot{\velocity}_{\vec{z}} : \nabla \velocity_{\vec{z}} + \nu \nabla \dot{\velocity}_{\vec{z}} : \nabla \adjvelocity_{\vec{z}} + ( (\dot{\velocity}_{\vec{z}}  \cdot \nabla) \velocity_{\vec{z}} ) \cdot \adjvelocity_{\vec{z}} + ( (\velocity_{\vec{z}}  \cdot \nabla) \dot{\velocity}_{\vec{z}} ) \cdot \adjvelocity_{\vec{z}}   \\
	&\hphantom{\int_{D_{\hat{\vu}}}\,} + \adjpressure_{\vec{z}} \Div{\dot{\velocity}_{\vec{z}}} - \dot{\pressure}_{\vec{z}} \Div{\adjvelocity_{\vec{z}}} \dx \\
	&+ \intD{ \nu \nabla \velocity_{\vec{z}} : \nabla \dot{\adjvelocity}_{\vec{z}} + ( (\velocity_{\vec{z}}  \cdot \nabla) \velocity_{\vec{z}} ) \cdot \dot{\adjvelocity}_{\vec{z}} - \pressure_{\vec{z}} \Div{\dot{\adjvelocity}_{\vec{z}}} - \vec{f} \cdot \dot{\adjvelocity}_{\vec{z}} + \dot{\adjpressure}_{\vec{z}} \Div{\velocity_{\vec{z}}} } \\
	&+ \int_{D_{\hat{\vu}}} - \nu \left( \nabla \velocity_{\vec{z}} \nabla \vec{W} \right) : \nabla \velocity_{\vec{z}} - \nu \left( \nabla \velocity_{\vec{z}} \nabla \vec{W} \right) : \nabla \adjvelocity_{\vec{z}} - \nu \nabla \velocity_{\vec{z}} : \left( \nabla \adjvelocity_{\vec{z}} \nabla \vec{W} \right) \\
	&\hphantom{+\int_{D_{\hat{\vu}}}\,} - ((\nabla \vec{W} \velocity_{\vec{z}} \cdot \nabla) \velocity_{\vec{z}}) \cdot \adjvelocity_{\vec{z}} + \pressure_{\vec{z}} {\nabla \adjvelocity_{\vec{z}}}^\top : \nabla \vec{W} - (\nabla \vec{f} \vec{W}) \cdot \adjvelocity_{\vec{z}} \\
	&\hphantom{+\int_{D_{\hat{\vu}}}\,} - \adjpressure_{\vec{z}} {\nabla \velocity_{\vec{z}}}^\top : \nabla \vec{W} \dx \\
	&+ \int_{D_{\hat{\vu}}} \Div{\vec{W}} \left( \frac{\nu}{2} \nabla \velocity_{\vec{z}} : \nabla \velocity_{\vec{z}} + \nu \nabla \velocity_{\vec{z}} : \nabla \adjvelocity_{\vec{z}} + ((\velocity_{\vec{z}}  \cdot \nabla) \velocity_{\vec{z}} ) \cdot \adjvelocity_{\vec{z}} - \pressure_{\vec{z}} \Div{\adjvelocity_{\vec{z}}}  \right.  \\
	&\hphantom{+\int_{D_{\hat{\vu}}}\Div{\vec{W}}\ \ } \left. \vphantom{\frac{\nu}{2}}- \vec{f} \cdot \adjvelocity_{\vec{z}}+ \adjpressure_{\vec{z}} \Div{\velocity_{\vec{z}}} \right) \dx \\
	&+ \mu \left( \left( \hat{\vec{h}}(\hat{\vu}) + \frac{\vec{\lambda}}{\mu} \right) - \max \left(0, \hat{\vec{h}}(\hat{\vu}) + \frac{\vec{\lambda}}{\mu}  \right) \right) \\
	& \qquad \cdot \begin{pmatrix}
		\left[\intui{\vec{W} \cdot \vec{n}} \right]_{i\in\{1,\dots,s\}} \\[3pt]
		\left[\frac{1}{\vol(\hat{u}_i)} \intui{ \left( \vec{x} - \bary(\hat{u}_i) \right) \left( \vec{W} \cdot \vec{n} \right) } \right]_{i\in\{1,\dots,s\}} \\[3pt]
		\left[\frac{1}{\vol(\hat{u}_i)} \intui{ \left( \bary(\hat{u}_i)- \vec{x} \right) \left( \vec{W} \cdot \vec{n} \right) } \right]_{i\in\{1,\dots,s\}}
	\end{pmatrix}.
\end{align*}
Reordering terms reveals the weak forms of the state and adjoint equation (\eqref{eq:weak-formulation} and \eqref{eqn:ModelStokesEquationsMomentumAdj}, respectively) for test functions $\dot{\velocity},\dot{\adjvelocity}\in V(D_{\hat{\vu}})$ and $\dot{\pressure},\dot{\adjpressure}\in L^2(D_{\hat{\vu}})$, which leads to
\begin{align}
\label{eq:ShapeDerivativeLA}
\begin{aligned}
&\Deriv \hat{L}_A (\hat{\vu}, \vec{\lambda}, \vec{z};\mu)\left[\vec{W}\right]	\\
&= \int_{D_{\hat{\vu}}} - \nu \left( \nabla \velocity_{\vec{z}} \nabla \vec{W} \right) : \left( \nabla \velocity_{\vec{z}} + \nabla \adjvelocity_{\vec{z}} \right) - \nu \left( \nabla \adjvelocity_{\vec{z}} \nabla \vec{W} \right) : \nabla \velocity_{\vec{z}} \\
&\hphantom{=\int_{D_{\vu}}\;} - ((\nabla \vec{W} \velocity_{\vec{z}} \cdot \nabla ) \velocity_{\vec{z}}) \cdot \adjvelocity_{\vec{z}} + \left( \pressure_{\vec{z}} {\nabla \adjvelocity_{\vec{z}}}^\top - \adjpressure_{\vec{z}} {\nabla \velocity_{\vec{z}}}^\top \right) : \nabla \vec{W}  - (\nabla \vec{f} \vec{W}) \cdot \adjvelocity_{\vec{z}} \\
	&\hphantom{=\int_{D_{\vu}}\;} + \Div{\vec{W}} \left( \nu \nabla \velocity_{\vec{z}} : \left( \frac{1}{2} \, \nabla \velocity_{\vec{z}} + \nabla \adjvelocity_{\vec{z}} \right) + ((\velocity_{\vec{z}}\cdot\nabla)\velocity_{\vec{z}})\cdot \adjvelocity_{\vec{z}}\right. \\
	&\hphantom{=\int_{D_{\vu}} + \Div{\vec{W}} \quad} \left. \vphantom{\frac{1}{2}}- \pressure_{\vec{z}} \Div{\adjvelocity_{\vec{z}}} -\vec{f} \cdot \adjvelocity_{\vec{z}} + \adjpressure_{\vec{z}} \Div{\velocity_{\vec{z}}} \right) \dx \\
	&\quad\ + \mu \left( \left( \hat{\vec{h}}(\hat{\vec{u}}) + \frac{\vec{\lambda}}{\mu} \right) - \max \left(0, \hat{\vec{h}}(\hat{\vec{u}}) + \frac{\vec{\lambda}}{\mu} \right) \right)\\
	&\qquad\quad\ \cdot   \begin{pmatrix}
	\left[\intui{\vec{W} \cdot \vec{n}} \right]_{i\in\{1,\dots,s\}} \\[3pt]
	\left[\frac{1}{\vol(\hat{u}_i)} \intui{ \left( \vec{x} - \bary(\hat{u}_i) \right) \left( \vec{W} \cdot \vec{n} \right) } \right]_{i\in\{1,\dots,s\}} \\[3pt]
	\left[\frac{1}{\vol(\hat{u}_i)} \intui{ \left( \bary(\hat{u}_i)- \vec{x} \right) \left( \vec{W} \cdot \vec{n} \right) } \right]_{i\in\{1,\dots,s\}}
\end{pmatrix}.
\end{aligned}
\end{align}

\section{Numerical investigations} \label{sec:numericalResults}
In Section~\ref{sec:AugLagMethod}, we first formulate the stochastic augmented Lagrangian method proposed in \cite{geiersbach2023stochastic} for problems on manifolds and briefly review its theoretical properties.  
Then, we provide algorithmic details for the numerical application in Section~\ref{sec:deterministicAlgorithm}. Section~\ref{sec:numerical-results} is dedicated to presenting numerical results.

\subsection{The stochastic augmented Lagrangian method on manifolds}
\label{sec:AugLagMethod}
\subsubsection{Formulation of the method} 
In this section, we describe the stochastic augmented Lagrangian method for solving problems of the form \eqref{eq:SO-problem-abstract-extended}. We collect the constraints in a vector $\vec{h}\colon\mathcal{M} \rightarrow \R^n$, $\vec{u} \mapsto \vec{h}(\vec{u}) = (h_1 (\vec{u}), \dots, h_n(\vec{u}))^\top$. 
As in \cite{geiersbach2023stochastic}, we define the parametrized Lagrangian and augmented Lagrangian by
\begin{align*}
& L(\vec{u}, \vec{\lambda}, \vec{z}) = J(\vec{u},\vec{z}) + \vec{\lambda}^\top \vec{h}(\vec{u}), \\
& L_A(\vec{u},\vec{\lambda}, \vec{z};\mu)\coloneqq J(\vec{u},\vec{z})+  \frac{\mu}{2} \Big \lVert \max \left( 0, \vec{h}(\vec{u})+ \frac{\vec{\lambda}}{\mu}\right) \Big \rVert_2^2 - \frac{\lVert \vec{\lambda}\rVert_2^2}{2\mu},
\end{align*}
respectively. Additionally, we define a feasibility measure and its induced sequence by
\begin{equation*}
\label{eq:definition-H}
H(\vec{u},\vec{\lambda}; \mu) \coloneqq \left\lVert \vec{h}(\vec{u}) - \max \left( 0, \vec{h}(\vec{u})+ \frac{\vec{\lambda}}{\mu}\right) \right\rVert_2, \quad H_{k}\coloneqq H(\vec{u}^{k},\vec{w}^{k-1};\mu_{k-1}).
\end{equation*}
In order to formulate a gradient-based optimization procedure on a Riemannian (product) manifold $(\mathcal{M}, \mathcal{G})$, we need to represent gradients with respect to the Riemannian metric $\mathcal{G}$ under consideration to define descent directions as well as the multi-exponential map to define the next (shape) iterate (cf., e.g., \cite{Geiersbach2022}). A brief introduction in these objects can be found, for example, in \cite[section 2.1]{geiersbach2023stochastic} or \cite{Geiersbach2022}. In the following, $\nabla_{\vec{u}} L_A$ denotes the gradient of the parametrized augmented Lagrangian with respect to $\mathcal{G}$. In addition, the exponential map at a point $\vec{u}\in \mathcal{M}$ is a mapping $\exp_{\vu}$ from the tangent space $T_{\vu}\mathcal{M}$ to the manifold $\mathcal{M}$.
We now have the components needed to formulate the proposed method, Algorithm~\ref{alg:algorithmAugLagr}, to solve \eqref{eq:SO-problem-abstract-extended}.

\begin{algorithm}
	\caption{Stochastic augmented Lagrangian method for shape optimization}
	\label{alg:algorithmAugLagr}
	\begin{algorithmic}[1]
		\State \textbf{Input:} Initial shape vector $\vu^1$, 
		 parameters $\gamma > 1$, $\tau \in (0,1)$, bounded set $B \subset \R^n$ 
		\State \textbf{Initialization:} $\mu_1 > 0$, $\vec{\lambda}^1 \in \R^n$,  $k\coloneqq1$ \setlength\itemsep{0ex}
		\While{$\vec{u}^k$, $\vec{\lambda}^k$ not converged}
		\State Choose $\vec{w}^k \in B$, step size $t_k$, iteration limit $N_k$, and batch size $m_k$
		\State ${\vec{u}^{k,1}}\coloneqq{\vec{u}^k}$
		\For{$j=1, \dots, N_k$}
		\State Generate i.i.d.~samples $\{\vec{\xi}^{k,j,1}, \dots, \vec{\xi}^{k,j,m_k}\} $ independent from prior samples
			\State $\vec{u}^{k,j+1}\coloneqq\exp_{\vec{u}^{k,j}}(-  \frac{t_k}{m_k}\sum_{l=1}^{m_k}\nabla_{\vec{u}} L_A(\vec{u}^{k,j},\vec{w}^k, \vec{\xi}^{k,j,l};\mu_k))$
			\EndFor
	\State ${\vec{u}^{k+1}}\coloneqq{\vec{u}^{k,j+1}}$
	\State  $\vec{\lambda}^{k+1}\coloneqq\mu_k \left(\vec{h}(\vec{u}^{k+1}) + \frac{\vec{w}^k}{\mu_k}  - \max\left(0,  \vec{h}(\vec{u}^{k+1}) +  \frac{\vec{w}^k}{\mu_k} \right) \right)$
		\IfThenElse{$H_{k+1} \leq \tau H_k$ or $k=1$ satisfied}{set $\mu_{k+1} = \mu_k$}{set $\mu_{k+1}\coloneqq\gamma \mu_k$}
		\State $k\coloneqq k+1$
		\EndWhile
	\end{algorithmic}
\end{algorithm}

Algorithm~\ref{alg:algorithmAugLagr} is based on the augmented Lagrangian method from \cite{Steck2018,Kanzow2018} and the stochastic gradient method from \cite[Algorithm~3]{Geiersbach2022}, combined with a minibatch stochastic gradient method. 
The sequence $\vec{w}^k$ is taken from a bounded set $B$; one choice would be $\vec{w}^k = \pi_B(\vec{\lambda}^k)$, with $\pi_B$ denoting the projection onto the box constraint set $B$. The algorithm's performance and complexity are highly dependent on the choices of step sizes $t_k$, iteration limits $N_k$, and batch sizes $m_k$.

\subsubsection{Choice of method parameters}
Here, we will make some formal arguments and refer to \cite{geiersbach2023stochastic} for further details. Let $f_k(\vec{u}) =  \E[L_A(\vec{u}, \vec{w}^k, \vec{\xi}; \mu_k)]$ and $\tilde{j}(\vu) = \E[J(\vec{u},\vec{\xi})]$ and suppose the random sequence of iterates $\{\vec{u}^k\}$ generated by Algorithm~\ref{alg:algorithmAugLagr} almost surely (a.s.) remain in a bounded set $\hat{B}$. Let $P_{1,0}\colon T_{\gamma(1)}\mathcal{U}^N \rightarrow T_{\gamma(0)}\mathcal{U}^N$ denote the parallel transport along the (unique) geodesic such that $\gamma(0) = \vec{u}$ and $\gamma(1) = \tilde{\vec{u}}$ and set $g_i({\vec{u}})\coloneqq h_i({\vec{u}}) + \frac{w^k_i}{\mu_k}-\max(0,h_i({\vec{u}}) + \frac{w^k_i}{\mu_k}).$ Let $\nabla \tilde{j}$ and $\nabla \vec{h}$ and be Lipschitz continuous; then, it is possible to obtain $L_k$-Lipschitz continuity of $\nabla f_k$ by splitting terms and observing
\begin{align*}
& \lVert P_{1,0} \nabla f_k(\tilde{\vu}) - \nabla f_k(\vu)\rVert_{\mathcal{G}^N}\\
&\leq \lVert  P_{1,0} \nabla \tilde{j}(\tilde{\vu}) -\nabla \tilde{j}(\vec{u}) \rVert_{\mathcal{G}^N}  + \mu_k \Big\lVert \sum_{i=1}^n P_{1,0}\nabla h_i(\tilde{\vu})  g_i(\tilde{\vu}) - \nabla h_i(\vu)  g_i(\vu) \Big\rVert_{\mathcal{G}^N}\\
&  \leq (L_{\tilde{j}} +\mu_k L_{\vec{h}}) \mathrm{d}(\tilde{\vu},\vu),
\end{align*} 
where $L_k\coloneqq L_{\tilde{j}} + \mu_k L_{\vec{h}}$ with $L_{\tilde{j}}$ and $L_{\vec{h}}$ being constants coming from $\tilde{j}$ and $\vec{h}$, respectively. By assuming that iterates are contained in $\hat{B}$, the Lipschitz continuity of $f_k$ only varies through the parameter $\mu_k$.
In \cite[Theorem 2.1]{geiersbach2023stochastic}, under various technical assumptions, asymptotic convergence of Algorithm~\ref{alg:algorithmAugLagr} was established using an adaptation of the randomly stopped method from \cite{Ghadimi2013,Ghadimi2014}. We showed that if the step size is chosen such that $t_k = \alpha_k/L_k$, $\alpha_k \in (0,2)$, we have
\begin{equation}
\label{eq:efficiency}
\frac{1}{N_k} \sum_{\ell=1}^{N_k}{\E}[\lVert \nabla f_k(\vec{u}^{k,\ell})\rVert_{\mathcal{G}}^2] \leq \frac{2L_k(f_k(\vec{u}^k)-f_k^*)}{(2\alpha_k-\alpha_k^2)N_k}+ \frac{\alpha_k M^2}{(2-\alpha_k)m_k} \quad \forall k,
\end{equation}
where $f_k^*\coloneqq\inf_{\vu \in \hat{B}} f_k(\vec{u})$, and $M$ is a constant satisfying $\E[\lVert \nabla_{\vec{u}} J(\vu,\vec{\xi}) - \nabla j(\vec{u})\rVert_{\mathcal{G}}^2] \leq M^2$ for all $\vu \in \mathcal{M}$. 
Additionally, if $L=\sup_k L_k<\infty$ and $N_k$, $\alpha_k$, as well as $m_k$ are chosen to satisfy
\begin{equation}
\label{eq:finite-sum}
\sum_{k=1}^\infty \frac{1}{(2\alpha_k-\alpha_k^2)N_k}+ \frac{\alpha_k}{(2-\alpha_k)m_k}  < \infty,
\end{equation}
then we have $\lVert  \nabla f_k(\vec{u}^{k+1})\rVert_{\mathcal{G}} \rightarrow 0$ a.s. The choice of $N_k$ and $m_k$ depend on the desired convergence rate for $\{\frac{1}{N_k} \sum_{j=1}^{N_k}{\E}[\lVert \nabla f_k(\vec{u}^{k,j})\rVert_{\mathcal{G}}^2]\}$. For instance, from \eqref{eq:efficiency}, it is clear that a choice of $m_k \propto N_k \propto k^\gamma$ ($\propto$ meaning up to a constant) yields a.s.~sublinear convergence for  $\gamma>1$; the choice $m_k \propto N_k \propto 2^k$ ensures~a.s.~linear convergence; and taking $m_k \propto N_k \propto k!$ results in a.s.~superlinear convergence.

\subsection{Algorithmic details regarding the inner loop} \label{sec:deterministicAlgorithm}
The shapes according to Algorithm \ref{alg:algorithmAugLagr} are updated by the exponential map. This computation is prohibitively expensive in most applications because a calculus of variations problem must be solved or the Christoffel symbols need be known. 
In the following, we consider $\mathcal{M}= M_s(\mathcal{U}^N)$ and remember that any element of $M_s(\mathcal{U}^N)$ can be understood as an element of $\mathcal{U}^N$.
An approximation using a (multi)-retraction
\[
\mathcal{R}_{\vec{u}^{k,j}}^N\colon T_{\vec{u}^{k,j}}\mathcal{U}^N\to \mathcal{U}^N,\, \vec{v}=(v_1,\dots,v_N)\mapsto (\mathcal{R}_{u^{k,j}_1}v_1,\dots,\mathcal{R}_{u^{k,j}_N}v_N)
\]
is often used
to update the shape vector $\vec{u}^{k,j}=(u^{k,j}_1,\dots,u^{k,j}_N)$.
If we consider $\mathcal{U}^N=B_e([0,1],\mathbb{R}^2)^N $, for each shape $u_i^{k,j}$, we can define the retraction in \cite[(3.14), page 265]{SchulzWelker} also on $B_e([0,1],\mathbb{R}^2)$ and use it in our computations.
In addition to updating a shape, we also need to update the computational mesh in each iteration.
In order to compute a deformation field that can be applied to $D_{\hat{\vu}}$ that  corresponds to the considered retraction, it has been described in~\cite[Eq.~(34)] {Geiersbach2022}
that the shape derivative can be used in an `all-at-once'-approach as the right-hand side of a variational problem, which is based on the Steklov--Poincar\'{e} metric\footnote{In this paper, we focus on the Steklov--Poincar\'{e} metric due to its advantages for the quality of the discretization of $D_{\hat{\vu}}$ (cf.~\cite{Schulz2016a,Siebenborn2017}).}: Find $\vec{V} \in \mathcal{W}(D_{\hat{\vu}})$ s.t.
\begin{equation}
	a_{\hat{\vu}}(\vec{V}, \vec{W}) = \Deriv \hat{L}_A (\hat{\vu}, \vec{\lambda},\vec{z};\mu)[\vec{W}] \quad \forall \vec{W}\in \mathcal{W}(D_{\hat{\vu}}),
	\label{deformatio_equation}
\end{equation}
where $a_{\hat{\vu}}$ is a coercive and symmetric bilinear form defined on $\mathcal{W}(D_{\hat{\vu}}) \times \mathcal{W}(D_{\hat{\vu}})$. The computed deformation field $\vec{V}$ from solving \eqref{deformatio_equation} is then applied to $D_{\hat{\vu}}$ as in~\eqref{eq:PerturbationOfIdentity}. 

For the bilinear form in \eqref{deformatio_equation}, we use 
\begin{equation}
\label{eq:BLF-elasticity}
a_{\hat{\vu}}^{\hat{\mu}}(\vec{V},\vec{W}) = \intD { 2 \hat{\mu} \varepsilon(\vec{V}):\varepsilon(\vec{W}) },
\end{equation}
in our simulations, where $\varepsilon(\vec{W})\coloneqq\tfrac{1}{2} (\nabla \vec{W} + \nabla \vec{W}^\top)$ and $\hat{\mu}$ is determined by solving
\begin{equation}
\label{eq:Poisson}
\Delta \hat{\mu} = 0 \quad \text{in } D_{\hat{\vu}},\qquad\qquad\hat{\mu} = \mu_{\max} \quad\text{on } \hat{\vu}, \qquad\qquad 
\hat{\mu}  = \mu_{\min} \quad \text{on } \Gamma,
\end{equation}
with the choices $\mu_{\max} = 33$ and $\mu_{\min} = 10$.

\subsection{Numerical results}
\label{sec:numerical-results}
All numerical simulations were performed either on a workstation equipped with  
$32$ physical cores, or the HPC cluster HSUper\footnote{Further information about the technical specifications can be found at \url{https://www.hsu-hh.de/hpc/en/hsuper/}.}.
In both cases Python 3.10.10 together with the FEniCS toolbox~\cite{AlnaesEtal2015}, version 2019.1.0, were used. The hold-all domain~$D$ was chosen to be $(-10,20) \times (-10,10)$ and the initial configuration and numbering of $s=5$ shapes can be seen in Figure~\ref{fig:shapes_iter0}. The two-dimensional computational mesh was generated using Gmsh 4.11.1 \cite{Geuzaine2009}. 
The resulting mesh contains 341 edges discretizing the boundaries, and 6611 triangle elements discretizing $D_{\hat{\vu}}$. The mesh was newly generated without requiring user interaction in case the mesh quality deteriorated.\footnote{We call the mesh quality deteriorated if the FEniCS function \texttt{MeshQuality.radius\_ratio\_min\_max} yields a value below $40\%$.} If the shape space definition is  based on the number of nodes belonging to the discretization of the shapes (see Remark~\ref{rem:KinksInShapes}) then the remeshing requires the definition of a new shape space if the number of nodes that discretize the shape is increased. To ensure inf-sup stability, stable P2-P1 Taylor-Hood elements were used for the discretization, see, e.g.,~\cite{Taylor1974,Bercovier1979}. Systems of equations were solved using the default sparse LU decomposition in combination with a Newton solver, as implemented in FEniCS. Random values were generated with \texttt{numpy} 1.22.4 using \texttt{numpy.random} with the seed $863860$, and parallelization of multiple stochastic realizations was performed using the \texttt{mpi4py} module, version 4.1.3 on the cluster and 4.1.4 on the workstation.

For all experiments, the parameters from Algorithm~\ref{alg:algorithmAugLagr} were chosen to be $\gamma = 2$ and $\tau = 0.9$ and penalty terms/multipliers were initialized to $\mu^1=1$, and $\vec{\lambda}^1 = \vec{0}$. Additionally, we projected the Lagrange multipliers $\vec{\lambda}^k$ onto $B=(-100,100)^{25}$ to obtain $\vec{w}^k$. The fluid viscosity was set to be $\nu=0.2$ and volumetric forces were neglected ($\vec{f}=\vec{0}$). The volume of each shape was constrained to be at or above $100\%$ initial volume, i.e. $\underline{\mathcal{V}}_i = \vol({\hat{u}_i^1})$ for every $i$. The barycenter locations were allowed to vary in the box $[-0.2, 0.5] \times[-0.3, 0.4]$ centered at the respective initial positions ${\vec{\mathcal{B}}_{1}^1}=(-0.5,5.5)^\top$, ${\vec{\mathcal{B}}_{2}^1}=(4.5, 0.5)^\top$, ${\vec{\mathcal{B}}_3^1}=(-5.5,0.5)^\top$,  ${\vec{\mathcal{B}}_{4}^1}=(-4.5, -5)^\top$, and ${\vec{\mathcal{B}}_{5}^1}=(2.5,-7)^\top$.

\subsubsection{The stochastic model}
For the velocity distribution $\vec{g}$ on the boundary~$\hat{\vu}$, we chose a no-slip condition $\vec{g}|_{\hat{\vec{u}}}=\vec{0}$. On~$\Gamma_D$, we defined a random field as in~\cite{Benner2020} by
\begin{equation}
\label{eq:def-g}
\vec{g}(\vx,\vec{\xi}(\omega)) = \begin{cases}
\begin{pmatrix}
\kappa(\vx,\vec{\xi}(\omega))\\
0
\end{pmatrix}, &\vx \in \{-10\}\times (-10,10),\\
\vec{0} & \text{else}
\end{cases}
\end{equation}
with 
\begin{equation}
\label{eq:def-kappa}
\kappa(\vec{x},\vec{\xi}(\omega)) = \frac{1}{100}(10+x_2)(10-x_2) + \sum_{\ell=1}^{20} \ell^{-\eta-1/2} \sin \left( \frac{\pi \ell x_2}{10}\right)\xi_\ell(\omega),
\end{equation}
where $\xi_\ell \sim U[-1,1]$ ($U[a,b]$ being the uniform distribution on the interval $[a,b]$)
and $\eta$ is a constant related to the smoothness of the (truncated) random field. The expected value of this random field is represented by the first term in \eqref{eq:def-kappa}, which describes a parabola-shaped inflow profile with a value of $1$ at the middle of the boundary and $0$ at the corners. The smoothness constant was chosen to be $\eta=2.5$. The $20 \times m_k$ matrix of random values was generated row by row; each column of this matrix, which represented a sample of the random field, was distributed to the $m_k$ processes.

\subsubsection{Estimating the Lipschitz constant $L_k$}
A crucial component in Algorithm~\ref{alg:algorithmAugLagr} is the step size $t_k$, which depends on the Lipschitz constant $L_k$. Since these $L_k$ (which depend on the gradients of the expectation $\hat{f}_k(\hat{\vec{u}})\coloneqq\E[\hat{L}_A(\hat{\vu}, \vec{w}^k, \vec{\xi};\mu_k)]$) are not available to us, we used an offline numerical procedure to estimate their values. We used a fixed number of samples $m=32$ corresponding to the number of cores available for parallel processing on the workstation. Given an outer iterate $k$ and an inner iterate $j$, we drew $m$ samples $\{ \vec{\xi}^{k,j,l}\}_{l=1}^{m}$ and computed an averaged deformation field $\bar{\vec{V}}^{k,j}\coloneqq\frac{1}{m} \sum_{l=1}^{m} \vec{V}^{k,j,l}$, where $\vec{V}^{k,j,l}$ solves 
\begin{equation}
\label{eq:deformation-equation-numerics}
a_{\hat{\vu}^{k,j}}(\vec{V}, \vec{W}) = \Deriv \hat{L}_A (\hat{\vu}^{k,j}, \vec{w}^k,\vec{\xi}^{k,j,l};\mu_k)[\vec{W}] \quad \forall \vec{W}\in \mathcal{W}(D_{\hat{\vu}^{k,j}})
\end{equation}
 and $D_{\hat{\vu}^{k,j}}$ is $D_{\hat{\vu}}$ corresponding to the shape vector $\hat{\vu}^{k,j}$ at the $j$th iteration in the $k$th outer loop; see Algorithm~\ref{alg:algorithmAugLagr}. 
Then, Armijo backtracking was performed with respect to the average deformation field, i.e., we found the step size $t^{k,j} = t_0 \beta^{j'}$ ($t_0>0$) with the smallest nonnegative integer $j'$  satisfying
\begin{equation}
\label{eq:backtracking-rule}
\begin{aligned}
&\frac{1}{m} \sum_{l=1}^{m} \hat{L}_A(\hat{\vu}^{k,j'},\vec{w}^{k}, \vec{\xi}^{k,j,l};\mu_k) \\
&\quad \leq \frac{1}{m} \sum_{l=1}^{m} \hat{L}_A(\hat{\vu}^{k,j},\vec{w}^k, \vec{\xi}^{k,j,l};\mu_k) - \sigma t^{k,j'} \lVert \bar{\vec{V}}^{k,j} \rVert_{H^1(D_{\hat{\vu}^{k,j}},\R^2)}^2,
\end{aligned}
\end{equation}
where $\hat{\vu}^{k,j} \in {M}_s^c$ is the starting shape and  $\hat{\vu}^{k,j'} = \hat{\vu}^{k,j} - t_0\beta^{j'} \bar{\vec{V}}^{k,j}|_{\hat{\vu}^{k,j}}$ is a shape that has been perturbed by $t_0\beta^{j'} \bar{\vec{V}}^{k,j}$. For the test, we used $\sigma=10^{-4}$, $\beta = 0.9$ and $t_0=8$. It is known (see \cite[Section~1.2.2]{Bertsekas1999}) that for a function with a Lipschitz gradient $L_k$, the backtracking procedure defined by \eqref{eq:backtracking-rule} will converge for all step sizes $t_k = \alpha_k/L_k$ with $\alpha_k \in \left[0, 2(1-\sigma)\right]$. We use this property to obtain an estimate for $L_k$ by taking the minimum accepted step size over all inner iterations~$j$, i.e., $L_k = \frac{2(1-\sigma)}{\min_j t^{k,j}}$. Since $L_k = L_{\tilde{j}} + \mu_k L_{\vec{h}}$ (see Section~\ref{sec:AugLagMethod}), 
we estimated $L_{\tilde{j}}$ and $L_{\vec{h}}$ by choosing different $\mu_k$. A least squares fit using the values in Table~\ref{tab:penalty_factor_step_sizes} yielded an estimate of $L_{\tilde{j}}=0.42215$ and $L_{\vec{h}}=0.36036$ with an $R^2$ value of $0.99858$. Using this information, the step size in Algorithm~\ref{alg:algorithmAugLagr} was set to $t_k = L_k^{-1}=(0.42215 + 0.36036\mu_k)^{-1}$ for the following stochastic numerical experiments, i.e., $\alpha_k =1$ for all $k$.

\begin{table}[tbp]
	\centering
	\caption{Minimal accepted step sizes for different penalty factors.}
	\begin{tabular}{l|rrrrrrr}
		$\mu$  & 1   & 2   & 4   & 8   & 16  & 32  & 64  \\
		\hline
		$\min_j t^{1,j}$ & 4.7239 & 2.2594 & 1.0807 & 0.57432 & 0.30522 & 0.16220 & 0.086202
	\end{tabular}
	\label{tab:penalty_factor_step_sizes}
\end{table}

\subsubsection{Stochastic solutions}

For the following experiment, we chose the step sizes using the estimated Lipschitz constants, the batch sizes $m_k= 2^{k-1}$, and the iteration numbers $N_k=2^{k+2}$. 
We approximated the objective functional using $\bar{j}(\hat{\vu}^{k, N_k}) = \frac{1}{m_k} \sum_{l=1}^{m_k} \hat{J}(\hat{\vu}^{k, N_k}, \vec{\xi}^{k,N_k,l})$ and the stationarity measure on the left-hand side of \eqref{eq:efficiency} using 
\[S^k \coloneqq \frac{1}{N_k} \sum_{j=1}^{N_k} \| \bar{\vec{V}}^{k,j}  \|_{H^1(D_{\hat{\vu}^{k,j}},\R^2)}^2 = \frac{1}{N_k} \sum_{j=1}^{N_k}\Big\lVert \frac{1}{m_k} \sum_{l=1}^{m_k} \vec{V}^{k,j,l}  \Big\rVert_{H^1(D_{\hat{\vu}^{k,j}},\R^2)}^2,\]
 where $\vec{V}^{k,j,l}$ solves
\eqref{eq:deformation-equation-numerics}. In Table~\ref{tab:NumericalResultsStochastic}, these values are displayed for each $k$ for a numerical experiment with a total of 11 outer iterations.
The results confirm the expected linear convergence rate of the algorithm with $m_k$ and $N_k$.
\begin{table}[htbp]
	\centering
	\caption{Numerical results for the number of inner iterations~$N_k$, batch size~$m_k$, objective function estimate $\bar{j}(\hat{\vu}^{k, N_k})$, and the stationarity measure~$S^k$, over outer iterations~$k$.}
	\begin{tabular}{r|rrrr}
		\multicolumn{1}{c|}{$k$} & \multicolumn{1}{c}{$N_k$} & \multicolumn{1}{c}{$m_k$} & \multicolumn{1}{c}{$\bar{j}(\hat{\vu}^{k, N_k})$} & \multicolumn{1}{c}{$S^k$} \\
		\hline
		1	&	8		&	1			&	11.7279	&	$4.7836 \cdot 10^{-2}$	\\
		2	&	16		&	2			&	9.01430	&	$1.2340 \cdot 10^{-2}$	\\
		3	&	32		&	4			&	8.79418	&	$5.2774 \cdot 10^{-3}$	\\
		4	&	64		&	8			&	9.20942 &	$1.8561 \cdot 10^{-3}$	\\
		5	&	128		&	16			&	8.62501 &	$1.4991 \cdot 10^{-3}$	\\
		6	&	256		&	32			&	8.75685	&	$7.6928 \cdot 10^{-4}$	\\
		7	&	512		&	64			&	8.41575	&	$1.7457 \cdot 10^{-4}$	\\
		8	&	1024	&	128		&	8.52603	&	$9.7468 \cdot 10^{-5}$	\\
		9	&	2048	&	256		&	8.50631	&	$3.4349 \cdot 10^{-5}$	\\
		10	&	4096	&	512		&	8.58534	&	$1.5611 \cdot 10^{-5}$	\\
		11	&	8192	&	1024	&	8.43387	&	$8.3566 \cdot 10^{-6}$
	\end{tabular}
	\label{tab:NumericalResultsStochastic}
\end{table}


\begin{figure}[htbp]%
	\centering%
	\setlength\figureheight{.47\textwidth}%
	\setlength\figurewidth{.47\textwidth}%
	\begin{subfigure}[t]{.49\textwidth}%
		\centering%
				\includegraphics{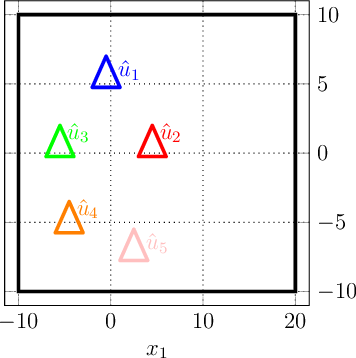}%
		\caption{Shapes before optimization. The $s=5$ shapes are numbered as: blue (1), red (2), green (3), orange (4) and pink (5).}%
		\label{fig:shapes_iter0}%
	\end{subfigure}%
	\hfill%
	\begin{subfigure}[t]{.49\textwidth}%
		\centering%
		\includegraphics{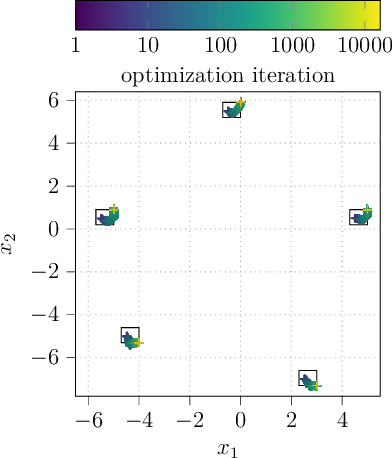}
		\caption{Position of the barycenter of each shape during stochastic optimization. The feasible barycenter positions are denoted by rectangles.}%
		\label{fig:NumericResults_stochastic_barycenter}%
	\end{subfigure}%
	\caption{Shapes after stochastic optimization (left) and barycenter positions over the course of the optimization (right).}
	\label{fig:NumericResults_stochastic}
\end{figure}

\begin{figure}[htbp]%
	\centering%
	\setlength\figureheight{.47\textwidth}%
	\setlength\figurewidth{.47\textwidth}%
	\begin{subfigure}[t]{.49\textwidth}%
			\centering%
					\includegraphics{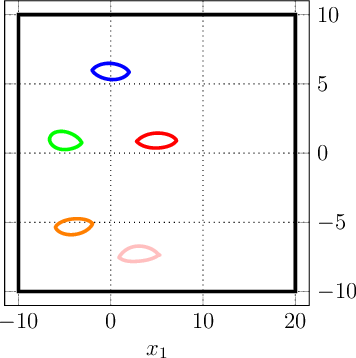}
			\caption{Shapes after stochastic optimization.}
			\label{fig:NumericResults_stochastic_shapesFinal}%
		\end{subfigure}\hfill%
	\begin{subfigure}[t]{.49\textwidth}%
			\centering%
					\includegraphics{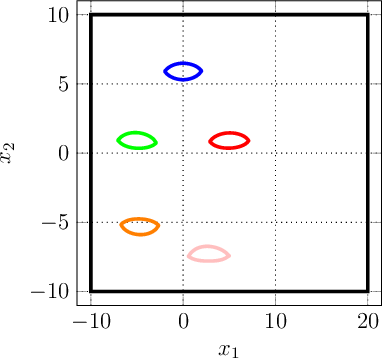}%
			\caption{Shapes after the deterministic optimization for $\xi_\ell=0$ $\forall \ell$.}%
			\label{fig:NumericResults_deterministic_shapesFinal}%
		\end{subfigure}
	\caption{Shapes $u_1$ (blue), $u_2$ (red), $u_3$ (green), $u_4$ (orange) and $u_5$ (pink) after the stochastic and deterministic optimization.}%
	\label{fig:shapes_initial_and_final}%
\end{figure}

The shapes at the start and end of the stochastic optimization are shown in Figure~\ref{fig:shapes_iter0} and Figure~\ref{fig:NumericResults_stochastic_shapesFinal}. All five shapes go from an initially triangular shape to a more streamlined shape. The initial kinks in the shape are removed by the optimization, and new kinks are formed towards the left and right of some shapes.

The barycenter position of each shape over the course of the optimization can be seen in Figure~\ref{fig:NumericResults_stochastic_barycenter}. One can see that all shapes move within the hold-all domain, and at some point of the optimization either the lower or upper bound is active for all shapes.

\begin{figure}[htb]
	\centering
	\setlength\figureheight{.25\textwidth}  
	\setlength\figurewidth{.25\textwidth}  
	\hfill\includegraphics[width=\textwidth]{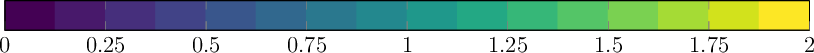}\\%
	\begin{subfigure}[b]{0.32\textwidth}%
		\includegraphics[scale=0.9]{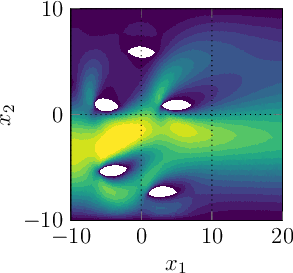}%
		\caption{$\hat{J}(\hat{\vu}^k, -\vec{1}) = 12.127$}%
		\label{fig:velocity_stochastic_Minus1}%
	\end{subfigure}\hspace*{.06\textwidth}%
	\begin{subfigure}[b]{0.32\textwidth}%
		\includegraphics[scale=0.9]{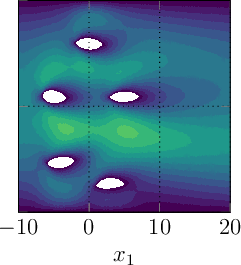}%
		\caption{$\hat{J}(\hat{\vu}^k, \vec{0}) = 6.8452$}%
		\label{fig:velocity_stochastic_Null}%
	\end{subfigure}
	\begin{subfigure}[b]{0.32\textwidth}%
		\includegraphics[scale=0.9]{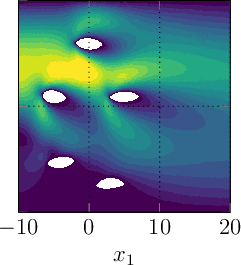}%
		\caption{$\hat{J}(\hat{\vu}^k, \vec{1}) = 11.518$}%
		\label{fig:velocity_stochastic_Plus1}%
	\end{subfigure}
	\caption{Fluid velocity magnitude with different inflow profiles \eqref{eq:def-g}--\eqref{eq:def-kappa} with the shapes from stochastic optimization (Figure~\ref{fig:NumericResults_stochastic_shapesFinal}). An estimate of the objective functional and the deformation field after optimization using random samples from seed 124764 with a sample size of $10016$ yields $\bar{j}(\hat{\vu}^{k})=8.4968$ and $S = \| \frac{1}{10016} \sum_{l=1}^{10016} \vec{V}^{l}  \|_{H^1(D_{\hat{\vu}^{k}},\R^2)}^2=3.2730 \cdot 10^{-6}$, where $\vec{V}^l$ solves \eqref{eq:deformation-equation-numerics}.}
\label{fig:velocity_stochastic}
\end{figure}

Further, the fluid velocity magnitude for the shapes in Figure~\ref{fig:NumericResults_stochastic_shapesFinal} and an evaluation of the objective functional for an inflow condition \eqref{eq:def-g}--\eqref{eq:def-kappa} using $\xi_\ell=-1$ $\forall \ell$, $\xi_\ell=0$ $\forall \ell$ and $\xi_\ell=1$ $\forall \ell$ is shown in Figure~\ref{fig:velocity_stochastic}, which indicates the strong influence of the randomness on the flow profile.

\subsubsection{Comparison to deterministic solutions}
In this experiment, we demonstrate that neglecting stochastic information leads to suboptimal shapes. We modified the velocity distribution defined in \eqref{eq:def-g} so that $\xi_\ell=0$ $\forall \ell$ is chosen in \eqref{eq:def-kappa}. For this, we used Armijo backtracking, which means replacing lines $6$--$8$ in Algorithm~\ref{alg:algorithmAugLagr} with the step size control in \cite[Algorithm~2]{Geiersbach2022}. On $D_{\hat{\vu}}$, at outer iteration $k$ and inner iteration $j$, this procedure translates to determining $t^{k,j} = t_0^j \beta^{j'}$ ($t_0>0$) with the smallest nonnegative integer $j'$ satisfying \eqref{eq:backtracking-rule} with $\vec{\xi}^{k,j,l}=\vec{0}$ and $m=1$, where the deterministic deformation vector field $\bar{\vec{V}}^{k,j}$ is the solution to \eqref{eq:deformation-equation-numerics} with $\vec{\xi}^{k,j,l}=\vec{0}$. The parameters here were chosen to be $\sigma=10^{-4}$, $\beta=\frac{1}{2}$, and $t_0^j=8$ (for $j \leq 10$). To reduce the number of rejected Armijo steps, we chose $t_0^j = \min_{\iota=1, \dots, 10} t_0^\iota$ for $j > 10$. For the inner loop termination condition, with the Lagrangian
 $\hat{L}(\hat{\vec{u}},\vec{\lambda}, \vec{z})\coloneqq\hat{J}(\hat{\vu},\vec{z}) + \vec{\lambda}^\top \hat{\vec{h}}(\hat{\vec{u}})$, we used the optimality measure $\hat{r}_k\coloneqq\hat{r}(\hat{\vec{u}}^k, \vec{\lambda}^k)$, where $\hat{r}(\hat{\vu},\vec{\lambda}) = \lVert \nabla_{\hat{\vec{u}}} \hat{L}(\hat{\vu},\vec{\lambda},\vec{0}) \rVert_{H^1(D_{\hat{\vec{u}}},\R^2)} + \lVert \hat{\vec{h}}(\hat{\vu}) + \max(0, \hat{\vec{h}}(\hat{\vu}) + \vec{\lambda})\rVert_2$.
The inner loop was terminated if $\frac{\hat{r}_{k+1}}{\hat{r}_k} \leq \frac{1}{k^2+2}$, a rule motivated by \cite[Theorem 2.3]{Steck2018}.

The shapes $\hat{u}_1, \ldots, \hat{u}_5$ before optimization are the same as in the previous experiment, and the shapes after 12 outer (augmented Lagrange) iterations are shown in Figure~\ref{fig:NumericResults_deterministic_shapesFinal}. It can be seen that the obtained barycenters varied significantly between the stochastic and deterministic optimization. This deterministic optimization also yielded a kink for shape $\hat{u}_3$ and $\hat{u}_4$ towards the left, while this area is more rounded after the stochastic optimization.

\begin{figure}[htb]
	\centering
	\setlength\figureheight{.25\textwidth}  
	\setlength\figurewidth{.25\textwidth}  
	\hfill\includegraphics[width=\textwidth]{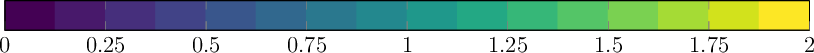}\\%
	\begin{subfigure}[t]{0.32\textwidth}%
		\includegraphics[scale=0.9]{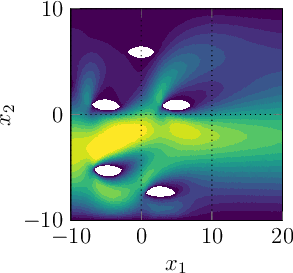}%
		\caption{$\hat{J}(\hat{\vu}^k, -\vec{1}) = 12.696$}%
		\label{fig:velocity_deterministic_Minus1}%
	\end{subfigure}\hspace*{.06\textwidth}%
	\begin{subfigure}[t]{0.32\textwidth}%
		\includegraphics[scale=0.9]{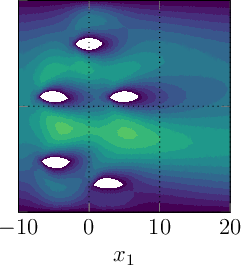}%
		\caption{$\hat{J}(\hat{\vu}^k, \vec{0}) = 6.7708$}%
		\label{fig:velocity_deterministic_Null}%
	\end{subfigure}
	\begin{subfigure}[t]{0.32\textwidth}%
		\includegraphics[scale=0.9]{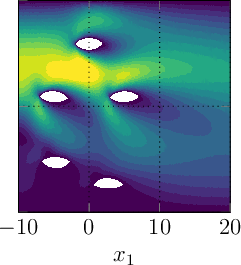}%
		\caption{$\hat{J}(\hat{\vu}^k, \vec{1}) = 11.818$}%
		\label{fig:velocity_deterministic_Plus1}%
	\end{subfigure}
	\caption{Fluid velocity magnitude with different inflow profiles \eqref{eq:def-g}--\eqref{eq:def-kappa} with the shapes from deterministic optimization (Figure~\ref{fig:NumericResults_deterministic_shapesFinal}). An evaluation of the objective functional and the deformation field with the same samples as in Figure~\ref{fig:velocity_stochastic} yields $\bar{j}(\hat{\vu}^{k})=8.5962$ and $S=7.1740 \cdot 10^{-4}$.}
	\label{fig:velocity_deterministic}
\end{figure}

Additionally, we show the fluid velocity magnitude after 12 outer (augmented Lagrange) iterations in Figure~\ref{fig:velocity_deterministic}. Here, an inflow profile with $\xi_\ell=-1$ $\forall \ell$, $\xi_\ell=0$ $\forall \ell$ and $\xi_\ell=1$ $\forall \ell$ was applied. The velocity profile looks qualitatively similar to the profiles from stochastic optimization (Figure~\ref{fig:velocity_stochastic}); however, higher objective functional values can be observed for values of $\xi_\ell$ further away from zero. An estimate of the objective functional $\bar{j}(\hat{\vu}^k)$ with a sample size of 10016 yields a reduction by $1.2\%$ when using a stochastic approach. The stationarity measure of the optimization is estimated to be strongly superior in the stochastic case by a factor of $100$, which indicates that the deterministic solution is not an optimal solution for the stochastic optimization.

\section{Conclusions}
\label{sec:conclusions}

In this paper, we connected the differential-geometric structure of the space of piecewise smooth shapes to shape calculus. We defined the multi-material derivative and the (stochastic) multi-shape derivative in relation to the new shape space. An application to a multi-shape optimization problem was considered, where the expected viscous energy dissipation should be minimized and the physical system is governed by the Navier--Stokes equations under uncertainty with additional geometrical inequality constraints. Using our proposed framework, we calculate the stochastic multi-shape derivative for our example.

In the numerical experiments, we confirmed the convergence rate of the stochastic augmented Lagrangian method previously introduced in \cite{geiersbach2023stochastic}. We discussed various aspects of the method, including the choice of parameters. For our tests, we estimated the sequence of Lipschitz constants for the penalized stochastic optimization problems numerically. These constants were used in a step size control for a computationally more expensive stochastic experiment. We  showed the optimized shapes from the stochastic model and compared them to those obtained using deterministic models, where the deterministic optimization with respect to $\xi_\ell=0$ for all $\ell$  yielded worse results than a stochastic approach. 

It is likely that the step size rule could be further optimized, for instance using an online (as opposed to offline) approach to estimate Lipschitz constants. Furthermore, a less conservative estimate of the Lipschitz constant could enable faster progress of the optimization. Our method shows promise for other applications, including three-dimensional and/or transient problems as well as related problems with other geometric constraints.

\section*{Acknowledgement(s)}

Computational resources (HPC-cluster HSUper) have been provided by the project hpc.bw, funded by dtec.bw -- Digitalization and Technology Research Center of the Bundeswehr. dtec.bw is funded by the European Union -- NextGenerationEU.

\section*{Disclosure statement}

The authors report there are no competing interests to declare.

\section*{Funding}

This work has been partly supported by the state of Hamburg (Germany) within the Landesforschungsförderung under project ``Simulation-Based Design Optimization of Dynamic Systems Under Uncertainties'' (SENSUS) with project number LFF-GK11. 

\bibliographystyle{abbrvnat}
\bibliography{references.bib}

\end{document}